\documentclass[12pt, abstracton]{scrartcl}

\usepackage[margin=2.5cm]{geometry}

\usepackage{pstricks, pst-tree, pst-node} 
\usepackage[latin1]{inputenc}
\usepackage{amssymb}
\usepackage{xfrac}     
\usepackage{lmodern}
\usepackage{stmaryrd}
\usepackage{graphicx}
\usepackage{amsmath}  
\usepackage{amsthm}
\usepackage{wrapfig} 
\usepackage{color}  
\usepackage{microtype} 
\usepackage{setspace}
\usepackage{amsfonts}
\usepackage{wasysym}
\usepackage{cancel}
\usepackage{verbatim}
\usepackage{dsfont}
\usepackage{paralist}

\usepackage[small,nohug,heads=vee]{diagrams}
\diagramstyle[labelstyle=\scriptstyle]

\DeclareMathOperator{\R}{\mathbb{R}}
\DeclareMathOperator{\Z}{\mathbb{Z}}
\DeclareMathOperator{\Q}{\mathbb{Q}}
\DeclareMathOperator{\C}{\mathbb{C}}

\DeclareMathOperator{\N}{\mathbb{N}}

\newcommand{\qf}[1]{\langle #1\rangle}

\newcommand{\laurent}[1]{(\!(#1)\!)}
\newcommand{\Pfister}[1]{\langle\!\langle #1\rangle\!\rangle}

\newtheoremstyle{plain2}
  {10pt}   
  {10pt}   
  {\itshape}  
  {0pt}       
  {\bfseries} 
  {}         
  {5pt plus 1pt minus 1pt} 
  {}          
  
 \newtheoremstyle{beweis}
  {10pt}   
  {10pt}   
  {\normalfont}  
  {0pt}       
  {\bfseries} 
  {:}         
  {5pt plus 1pt minus 1pt} 
  {}          

\newtheoremstyle{definition2}
  {10pt}   
  {10pt}   
  {\normalfont}  
  {0pt}       
  {\bfseries} 
  {}         
  {5pt plus 1pt minus 1pt} 
  {}          

\theoremstyle{plain2}
\newtheorem{Satz}{Satz}[section]
\newtheorem{Lemma}[Satz]{Lemma}
\newtheorem{Proposition}[Satz]{Proposition}
\newtheorem{Theorem}[Satz]{Theorem}
\newtheorem{Korollar}[Satz]{Corollary}
\newtheorem{KorDef}[Satz]{Corollary and Definition}

\theoremstyle{definition2}
\newtheorem{Definition}[Satz]{Definition}
\newtheorem{Bemerkung}[Satz]{Remark}
\newtheorem{Beispiel}[Satz]{Example}

\theoremstyle{beweis}
\newtheorem*{Beweis}{Proof}

\begin{document}

%
%
%

\begin{titlepage}
	\vspace*{\fill}
	\centering
	{\scshape\LARGE\bfseries Pfister Numbers over Rigid Fields \par}
	\vspace{1.5cm}
	{\scshape Nico Lorenz\par}
	\vspace{1.5cm}
	{Fakult\"at f\"ur Mathematik, Technische Universit\"at Dortmund, D-44221 Dortmund, Germany\par}
	\vspace{2cm}
	{e-mail: nico.lorenz@tu-dortmund.de\par}
	\vspace{2cm}
	{\today\par}
	\vspace*{\fill}
\end{titlepage}

\begin{abstract}
	For certain types of quadratic forms lying in the $n$-th power of the fundamental ideal, we compute upper bounds and where possible exact values for the minimal number of general $n$-fold Pfister forms, that are needed to write the Witt class of that given form as the sum of the Witt classes of those $n$-fold Pfister forms. We restrict ourselves mostly to the case of so called rigid fields, i.e. fields in which binary anisotropic forms represent at most 2 square classes.
\end{abstract}


\begin{section}{Introduction}

Throughout this paper, let $F$ be a field of characteristic different from 2. By a quadratic form or just form for short, we will always mean a finite dimensional non-degenerate quadratic form over $F$. We will denote isometry of two forms $\varphi_1,\varphi_2$ by $\varphi_1\cong\varphi_2$. In abuse of notation, we will denote the Witt class of a quadratic form $\varphi$ again by $\varphi$. An $n$\textit{-fold Pfister form} for some $n\in\N$ is a form of the shape $\Pfister{a_1,\ldots, a_n}:=\qf{1, -a_1}\otimes\qf\ldots\otimes\qf{1,-a_n}$ with $a_1,\ldots, a_n\in F^\ast$. The set of $n$-fold Pfister forms over $F$ is denoted by $P_nF$, the set of forms that are similar to some $n$-fold Pfister form is denoted by $GP_nF$. By abuse of notation, we will call any form in $GP_nF$ a Pfister form. Both $P_nF$ and $GP_nF$ generate the $n$-th power of the \textit{fundamental ideal} $IF$, which we denote by $I^nF$, both as an additive group and as an ideal, i.e. any Witt class $\varphi\in I^nF$ can be written as a sum $\varphi=\pi_1+\ldots+\pi_m$ for suitable $\pi_1,\ldots, \pi_m\in GP_nF$. A frequently asked question in the current research on quadratic forms and the main topic of this article is now to determine the lowest number of Pfister forms that are needed to write $\varphi\in I^nF$ as such a sum. This number is called the $n$-Pfister number of $\varphi$. We will collect some known facts about it in the next section.\\
During this paper we will mostly consider rigid fields. Such fields admit complete discrete valuations, so section 3 will deal with some basic questions concerning Pfister numbers of forms over complete discrete valuation fields.\\
In section 4, we will then introduce rigid fields and study the behaviour of quadratic forms over such fields in a more detailed way. \\
Section 5 is devoted to both 14-dimensional forms in $I^3$ and 8-dimensional forms in $I^2$ over rigid fields. There is a strong connection between forms of these two types and the classification of them over rigid fields will help us to investigate 16-dimensional $I^3$-forms in the following section. We discuss these forms in detail as this is the smallest type of forms whose Pfister number is no upper bound available for so far.\\
The final part of this article now deals with the growth of Pfister numbers for increasing dimension.

\end{section}

\begin{section}{Basic Results on Pfister Numbers}

As already mentioned in the introduction, we would like to study the so called Pfister number, which we will now formally introduce.

\begin{Definition}
	We define the $n$\textit{-Pfister number} of a quadratic form $\varphi\in I^nF$ to be
	$$GP_n(\varphi):=\min\{k\in\N\mid \text{there are }\pi_1,\ldots, \pi_k\in GP_nF\text{ with }\varphi=\pi_1+\ldots+\pi_k\in WF\}.$$
	For a subset $S\subseteq WF$ and an integer $d\in\N$, we define
	$$GP_n(S, d):=\sup\{GP_n(\varphi)\mid\varphi\in S\cap I^nF, \dim\varphi\leq d\}.$$
	Additionally, we define the shortcuts
	$$GP_n(F, d):=GP_n(WF, d)~~\text{ and }~~GP_n(S):=\bigcup_{d\in\N}GP_n(S,d).$$
	We further define the \textit{unscaled} $n$\textit{-Pfister number} of $\varphi$ to be
	\begin{align*}
		{P}_n(\varphi):=\min\{k\in\N\mid &\text{ there are }\varepsilon_1,\ldots, \varepsilon_k\in\{\pm1\}\text{ and } \\
		&\pi_1,\ldots, \pi_k\in P_nF\text{ with }\varphi=\varepsilon_1\pi_1+\ldots+\varepsilon_k\pi_k\in WF\}.
	\end{align*}
	If the integer $n$ is clear from the context, we will often just say \textit{(unscaled) Pfister number}.
\end{Definition}

The main task in this article is now to calculate Pfister numbers in terms of invariants of a given form. As this seems to be a quite tough task, we will often be satisfied with upper or lower bounds. We will concentrate on the scaled version, as we have the following correspondence between both versions.

\begin{Proposition}\label{PropVergleichScaledUnscaled}
	For any quadratic form $\varphi$ over $F$ and any $n\in\N$, we have ${P}_n(\varphi)\leq2\cdot{GP}_n(\varphi).$
\end{Proposition}
\begin{Beweis}
	For any $a,x_1,\ldots, x_n\in F^\ast$, we have
	\begin{align*}
		a\Pfister{x_1,\ldots, x_n}&=\Pfister{x_1,\ldots, x_{n-1}}\otimes(a\Pfister{x_n})\\
		&=\Pfister{x_1,\ldots, x_{n-1}}\otimes(\qf{1, a}\perp -\qf{1,ax_n})\\
		&=\Pfister{x_1,\ldots, x_{n-1}, -a}\perp-\Pfister{x_1,\ldots, x_{n-1}, -ax_n}
	\end{align*}
	which then readily implies the assertion.
\end{Beweis}

For Pfister numbers in $I^2$, we have the following two results.

\begin{Proposition}{\cite[Chapter X. Exercise 4]{Lam2005}}\label{2PfisterZahl}
	Let $\varphi\in I^2F$ be a form of dimension $\dim\varphi\in\N$. Then $\varphi$ is Witt equivalent to a sum of at most $\frac{\dim\varphi}2-1$ forms in $GP_2F$.
\end{Proposition}

\begin{Beispiel}\label{Maximal2PfisterNumber}{(Parimala, Suresh, Tignol)}
	Let $K$ be a field and $F:=K\laurent{X_1}\ldots \laurent{X_n}$ for some $n\in\N$ with $n\geq 2$.
	According to (the proof of) \cite[Theorem 2.2]{ParimalaSureshTignol} (in which the assumption that $-1$ is a square is only needed to assure that the upcoming forms lie in $I^2F$ and can be omitted by adding a sign as below), we see that we have
	$${P}_2F(\qf{1, X_1,\ldots, X_n, (-1)^{\frac{n+2}2}X_1\cdot\ldots\cdot X_n})=n-1\text{ if }n\text{ is even}$$
	and
	$${P}_2F(\qf{X_1,\ldots, X_n, (-1)^{\frac{n+1}2}X_1\cdot\ldots\cdot X_n})=n-1\text{ if }n\text{ is odd}.$$
	Thus, \ref{PropVergleichScaledUnscaled} implies
	$${GP}_2F(\qf{1, X_1,\ldots, X_n, (-1)^{\frac{n+2}2}X_1\cdot\ldots\cdot X_n})\geq\frac{n-1}2$$
	and
	\begin{align}\label{LauraIstGerade}
		{GP}_2F(\qf{X_1,\ldots, X_n, (-1)^{\frac{n+1}2}X_1\cdot\ldots\cdot X_n})\geq\frac{n-1}2.
	\end{align}
	In the case where $n$ is even, using that the Pfister number is always an integer, we even get 
	\begin{align}\label{HierImZimmer}
		{GP}_2F(\qf{1, X_1,\ldots, X_n, (-1)^{\frac{n+2}2}X_1\cdot\ldots\cdot X_n})\geq\frac{n}2.
	\end{align}
	The reverse inequalities are covered in \ref{2PfisterZahl}, so we have equalities both in \eqref{LauraIstGerade} and \eqref{HierImZimmer}. Of course, since the values of ${GP}_2F$ are invariant under scaling and since we can redefine the indeterminates, we can restrict ourselves to the case where $n$ is even and just consider the form
	$$\varphi:=\qf{1, X_1,\ldots, X_n, (-1)^{\frac{n+2}2}X_1\cdot\ldots\cdot X_n}\in I^2F$$
	with $\dim\varphi=n+2$ and
	$${GP}_2(\varphi)=\frac n2,$$
	which is the biggest possible value. This form will also be referred to as the \textit{generic (rigid)} $I^2$-form of dimension $n+2$.
\end{Beispiel}

For further results concerning Pfister numbers, we refer the reader to \cite{HoffmannTignol}, \cite{IzhboldinKarpenko}, \cite{Karpenko2017} and \cite{BrosnanReichsteinVistoli}.

\end{section}

\begin{section}{Valuation Theoretic Results}

As the yet to be defined rigid fields, which we want to study, admit valuations, we want to give a short exposition of the main ingredients coming from valuation theory that will be used in the sequel. In this section, we fix a field $F$ equipped with a complete discrete valuation $v$ with residue class field $K$. We will often use the case of a Laurent series extension $F=K\laurent t$ with the usual valuation given by the least index such that the respective coefficient of the Laurent series is not zero.

\begin{Korollar}\label{PfisterDarstellungenLiften}
	Let $\varphi\in I^nF$ be a unimodular form. Then, the $n$-Pfister number of $\varphi$ over $F$ and of its first residue class form $\overline\varphi$ over $K$ coincide.
\end{Korollar}
\begin{Beweis}
	If we have $\overline \varphi=\overline\pi_1+\ldots+\overline\pi_k$ for some Pfister forms $\overline\pi_1,\ldots,\overline\pi_k\in GP_nK$, we can lift them to get a representation $\varphi = \pi_1+\ldots +\pi_k$ by \ref{LemmaSplitExactSequence}.\\
	For the converse, we fix a uniformizing element $t$. Using the isometry $\Pfister{at, bt}\cong\Pfister{at, ab}$ for all $a,b\in F^\ast$, we can find a representation 
	$$\varphi = \pi_1+\ldots +\pi_k+ \tilde\pi_1\otimes\Pfister{c_1t}+\ldots+\tilde\pi_\ell\otimes\Pfister{c_\ell t}+t\hat\pi_1+\ldots+t\hat\pi_m$$
	with unimodular forms $\pi_1,\ldots, \pi_k, \hat\pi_1,\ldots,\hat\pi_m\in GP_nF$ and $\pi_1,\ldots,\pi_\ell\in GP_{n-1}F$ and $c_1,\ldots, c_\ell\in F^\ast$.
	By comparing both residue class forms, we see that in $WF$, we have equalities
	$$\varphi = \pi_1+\ldots+\pi_k+c_1\tilde\pi_1+\ldots+c_\ell\tilde\pi_\ell ~~\text{and}~~ c_1\tilde\pi_1+\ldots+c_\ell\tilde\pi_\ell = \hat\pi_1+\ldots+\pi_m.$$
	This implies 
	$$\varphi = \pi_1+\ldots+\pi_k + \hat\pi_1+\ldots+\hat\pi_m,$$
	where all forms are unimodular. Thus, the claim follows.
\end{Beweis}

\begin{Proposition}\label{PropPfisterZahlUnimodularMalUniformizer}
	Let $\psi\in I^{n-1}K$ be a unimodular form and $\varphi:=\Pfister t\otimes\psi$ for some uniformizer $t$. We then have ${GP_{n-1}(\psi)=GP_n(\varphi)}$.
\end{Proposition}
\begin{Beweis}
	The inequality $GP_{n-1}(\psi)\geq GP_n(\varphi)$ is clear. For the converse, we consider a representation
	$$\varphi = \pi_1+\ldots +\pi_k+ \tilde\pi_1\otimes\Pfister{c_1t}+\ldots+\tilde\pi_\ell\otimes\Pfister{c_\ell t}+t\hat\pi_1+\ldots+t\hat\pi_m$$
	as above. After comparing residue class forms, we see that we have
	$$\pi_1+\ldots+\pi_k+\tilde\pi_1+\ldots+\tilde\pi_\ell=\psi= -\hat\pi_1-\ldots-\hat\pi_m+c_1\tilde\pi_1+\ldots+c_\ell\tilde\pi_\ell.$$
	These are representations of $\psi$ as a sum of $2k+\ell$ respectively $2m+\ell$ forms in $GP_{n-1}F$. If we had $k+\ell+m<GP_{n-1}F(\psi)$ one of the terms $2k+\ell$ and $2m+\ell$ would also be strictly smaller than $GP_{n-1}F(\psi)$, a contradiction. Thus, we have $GP_{n-1}(\psi)\leq GP_n(\varphi)$ and the proof is complete.
\end{Beweis}

\begin{Proposition}\label{DefinedOverFThenInPreimage}
	Let $\varphi$ be a quadratic form that lies in $I^nF(X)$ or $I^nF\laurent t$ defined over $F$. Then, there is a unique preimage $\psi\in WF$ under the canonical map $r_{F(X)/F}$ respectively $r_{F\laurent t/F}$ and it fulfills $\psi\in I^nF$.
\end{Proposition}
\begin{Beweis}
	We will denote the map induced by scalar extension in both cases by $r$.
	The existence and uniqueness of some $\psi\in WF$ with $r(\psi)=\varphi$ is clear as $\varphi$ is defined over $F$ and $r$ is known to be injective, see e.g. \cite[Chapter IX. Lemma 1.1]{Lam2005}.\\
	As $\varphi$ has a preimage in $I^nF$ because of \cite[Theorem 21.1, Corollary 21.3]{ElmanKarpenkoMerkurjev2008} respectively \cite[Exercise 19.15]{ElmanKarpenkoMerkurjev2008}, the claim follows.
\end{Beweis}

\begin{Korollar}\label{UeberLaurentBleibtGleich}
	Let $\varphi\in I^nF$ and $E$ be a field with $F(t)\subseteq E\subseteq F\laurent t$. We then have ${GP_n(\varphi)=GP_n(\varphi_E)}$.
\end{Korollar}
\begin{Beweis}
	As the Pfister number can only decrease when going up to a field extension, it is enough to show the inequality
	$$GP_n(\varphi)\leq GP_n(\varphi_E)$$
	for $E=F\laurent t$, but this follows directly from \ref{PfisterDarstellungenLiften}.
\end{Beweis}

In the sequel, we will frequently use the following well known exact sequence.

\begin{Lemma}{\cite[Exercise 19.15]{ElmanKarpenkoMerkurjev2008}}\label{LemmaSplitExactSequence}
	For all $n\in\N$ we have a split exact sequence
	\begin{align*}
		\begin{diagram}
			0 & \rTo & I^nK & \rTo & I^n F & \rTo & I^{n-1}K & \rTo & 0
		\end{diagram}
	\end{align*}
	where the maps are given by lifting and taking the second residue class form.
\end{Lemma}

The following result should be compared with \cite[Lemma 1.5]{Raczek2013}.

\begin{Proposition}\label{UnterteilenInUnimodulareFormen}
	Let $\varphi\in I^nF$ be a quadratic form such that both residue class forms are not hyperbolic. Then there is uniformizer $t$, unimodular forms $\sigma\in I^nF$ and $\tau\in I^{n-1}F$ with $\varphi=\sigma\perp\Pfister {-t}\otimes\tau\in WF$ and $\dim\sigma<\dim\varphi$.
\end{Proposition}
\begin{Beweis}
	We denote the first respectively second residue class forms of $\varphi$ with respect to some uniformizing element $t$ with $\varphi_1$ respectively $\varphi_2$. We then have
	\begin{align}\label{IchSagJBIhrSagtO}
		\varphi=\varphi_1\perp t\varphi_2=\varphi_1\perp-\varphi_2\perp\varphi_2\perp t\varphi_2=\varphi_1\perp -\varphi_2\perp\Pfister{-t}\otimes\varphi_2.
	\end{align}
	After multiplying $t$ with some unit of the valuation ring, i.e. changing the uniformizer, we can assume ${D_F(\varphi_1)\cap D_F(\varphi_2)\neq\emptyset}$. Then the form $\varphi_1\perp -\varphi_2$ is isotropic. If we choose 
	\begin{align*}
		\sigma:=(\varphi_1\perp-\varphi_2)_{\text{an}} \text{ and }
		\tau:=\varphi_2,
	\end{align*}
	we have $\dim\sigma<\dim\varphi$ and $\tau\in I^{n-1}F$ by \ref{LemmaSplitExactSequence}. Finally \eqref{IchSagJBIhrSagtO} implies ${\varphi\equiv \varphi_1\perp-\varphi_2\mod I^nF}$, which then leads to $\sigma:=(\varphi_1\perp-\varphi_2)_{\text{an}}\in I^nF$.
\end{Beweis}

With the above result, we are now in the position to bound the Pfister numbers of forms over a complete discrete valuation field in terms of Pfister numbers over the associated residue class field. As a first step, we record the following special case which follows directly by \ref{UnterteilenInUnimodulareFormen}.

\begin{Korollar}\label{PfisterzahlUnimodulareZerlegung}
	Let $\varphi$ be as in \ref{UnterteilenInUnimodulareFormen}. Then its $n$-Pfister number is bounded by 
	$${GP}_n\left(K, \dim(\varphi)-2\right)+{GP}_{n-1}\left(K, \frac12\dim(\varphi)\right).$$
\end{Korollar}
\begin{Beweis}
	We use the notation as in the proof of \ref{UnterteilenInUnimodulareFormen}. Since the Pfister number of any form is invariant under scaling, we can assume $\dim\varphi_2\leq \frac12\dim\varphi$. We thus get $\sigma\in I^nF$ and $\tau\in I^{n-1}F$ such that we have a representation $\varphi=\sigma +\Pfister{-t}\otimes \tau$ in the Witt ring $WF$ with some suitable uniformizer $t$ and 
	$$\dim\sigma\leq\dim\varphi-2 \text{ and } \dim\tau\leq \frac12\dim\varphi,$$
	where the first inequality can be assumed by \ref{UnterteilenInUnimodulareFormen} since both residue forms are not hyperbolic.
	By \ref{PfisterDarstellungenLiften} we have 
	$${GP}_n\left(\Pfister{-t}\otimes\tau\right)\leq {GP}_{n-1}\left(K, \frac12\dim(\varphi)\right)$$
	and the result now follows.
\end{Beweis}

As the main result of this section, we have the following:

\begin{Theorem}
	Let $F$ be complete discrete valuation field such that the characteristic of the residue class field $K$ is not equal to $2$. Then for all $n\in\N$ and all $d\in2\N$, we have
	\begin{align*}
		{GP}_n(F, d)\leq \max\left\{{GP}_n(K, d-2)+{GP}_{n-1}\left(K, \frac d2\right), {GP}_n(K, d)	\right\}.
	\end{align*}
\end{Theorem}
\begin{Beweis}
	For any $d$-dimensional quadratic form $\varphi\in I^nF$, either both of its residue class forms are not hyperbolic or $\varphi$ is similar to an unimodular form. The claim now follows by \ref{PfisterzahlUnimodulareZerlegung} and \ref{PfisterDarstellungenLiften}.
\end{Beweis}

We conclude this section with a remark in which we show how we can treat the case of a residue class field that admits a complete discrete valuation itself.

\begin{Bemerkung}\label{RigidFieldsTauscheBewertungen}
	Let $K$ be a field of characteristic not 2, $F=K\laurent{t_1}\laurent{t_2}$ and $\varphi$ a quadratic form over $F$. We then have the choice if we want to consider $\varphi$ as a form over $F$ or as a form over $E:=K\laurent{t_2}\laurent{t_1}$ as the $\mathbb F_2$-linear map $\Phi:F^\ast/F^{\ast2}\to E^\ast/E^{\ast2}$ defined by
	\begin{align*}
		aF^{\ast2}&\mapsto aE^{\ast2}~~~\text{for all }a\in K^\ast;\\
		t_1F^{\ast2}&\mapsto t_1E^{\ast2};\\
		t_2F^{\ast2}&\mapsto t_2E^{\ast2}
	\end{align*}
	is a group isomorphism with 
	\begin{align}
		\Phi(-1)=-1\text{ and }\Phi(D_F(\qf{x_1,\ldots, x_n}))=D_E\big(\qf{\Phi(x_1),\ldots, \Phi(x_n)}\big)\label{WaldSpaziergang}
	\end{align}
	for all $n\in\N$ and $x_1,\ldots, x_n\in F^\ast$, see \cite[Chapter XII. Harrison-Cordes Theorem 1.8]{Lam2005}. 
	It is an isomorphism as if $\{a_i\mid i\in I\}$ is a system of representatives of $K^\ast/K^{\ast2}$, then $\{a_i, a_it_1, a_it_2, a_it_1t_2\mid i\in I\}$ is a system of representatives of both $F^\ast/F^{\ast2}$ and $E^\ast/E^{\ast2}$.
	We would like to emphasize the fact that we have $F\neq E$ and that $\Phi$ is not an identity map even though it looks like one, especially when abusing the notation and identifying a non zero element of the field with its square class.
	\\
	We now further assume that we have $\varphi\in I^nF$ for some $n\in\N$. As $\Phi$ induces a ring isomorphism that obviously takes 1-fold Pfister forms over $F$ to 1-fold Pfister forms over $E$, the Pfister number of $\varphi$ over $E$ is lower than or equal to the Pfister number of $\varphi$ over $F$. As we can argue the other way round with $\Phi^{-1}$ we can even say that both Pfister numbers coincide.\\
	As the symmetric group $S_n$ for $n\geq2$ is generated by transpositions of the form $(k~k+1)$ for all $k\in\{1,\ldots, n-1\}$, we can extend the above to handle the field $K\laurent{t_1}\cdots\laurent{t_n}$ by reordering the Laurent variables in an appropriate way.\\
	As a last trick, we would like to mention that we can always change the uniformizing element in some ways. For example, we have $K\laurent t=K\laurent{at}$ for all $a\in K^\ast$, so for example $K\laurent{t_1}\laurent{t_2}=K\laurent{t_1}\laurent{t_1t_2}$.\\
	We will use these facts without mentioning them explicitly several times in the sequel. The main idea while using this is that quadratic forms are good to manage if they have a well understood subform. It is thus convenient to reorder the Laurent variables such that one gets a residue form of small dimension.
\end{Bemerkung}

\end{section}

\section{Introduction to the Theory of Rigid Fields}

Inspired by the work of M. Raczek \cite{Raczek2013}, we will prove upper bounds for the Pfister number of so called rigid fields. Using similar arguments, we generalize a lot of the arguments used in the just cited article. In the theory of quadratic forms, rigid fields are of interest because of several reasons. Firstly, they are simple enough to handle to build up a theory that already started in the late 1970s, see \cite{Ware}. As an example, there are a lot of interesting Galois-theoretic results available for rigid fields. Furthermore, nonreal rigid fields with a finite number of square classes are examples of the so called $\overline C$-fields. These are extreme examples as these are those fields that have the maximal number of anisotropic quadratic forms that can occur, when considering nonreal fields with finitely many square classes, see \cite[Chapter XI., Theorem 7.10, 7.14, Definition 7.16]{Lam2005}.

\begin{Definition}
	A field $F$ is called \textit{rigid}, if, for any binary anisotropic quadratic form $\beta$ over $F$, we have $|D_F(\beta)|\leq 2$.
\end{Definition}

\begin{Beispiel}
	As the square class groups of finite fields or euclidean fields consist of only two elements, these fields are rigid. Over a quadratically closed field there are no binary anisotropic forms. Thus quadratically closed fields are rigid as well.
\end{Beispiel}

We will now give a characterization of rigid fields that will be useful in the sequel.

\begin{Theorem}{\cite[Theorems 1.5, 1.8, 1.9]{Ware}}\label{RigidCharakterisierung}
	For a field $F$ the following are equivalent:
	\begin{enumerate}[(i)]
		\item $F$ is rigid;
		\item\label{EsIstMirEgaaaaaaaal} we have an isomorphism $WF\cong (\Z/n\Z)[G]$ with $n\in\{0, 2, 4\}$ and $G$ a group of exponent 2;
		\item\label{DochEsIstMirNichtsEingefaaaaaaallen} we have an isomorphism $WF\cong (\Z/n\Z)[H]$ with either $n=2$ and $H=F^\ast/F^{\ast2}$ or $n\in\{0, 4\}$ and $H\subseteq F^\ast/F^{\ast2}$ a subgroup with $-1\notin H$ and $[F^\ast/F^{\ast2}:H]=2$;
		\item for any anisotropic form $\varphi$, we have $|D_F(\varphi)|\leq \dim\varphi$;
		\item for any quadratic field extension $K/F$, the image of the inclusion map ${\iota:F^\ast/F^{\ast2}\to K^\ast/K^{\ast2}}$ has index $\leq2$. 
	\end{enumerate}
\end{Theorem}

An important field invariant when studying quadratic forms is the so called \textit{level} of a field, in symbols $s(F)$. It is defined as the least number $n$ of squares such that -1 is a sum of $n$ squares or $\infty$ if no such integer exists or equivalently the least integer $n$ such that $(n+1)\times\qf{1}$ is isotropic. It is well known that the level is either $\infty$ or a power of 2, see \cite[Chapter XI. Pfister's Level Theorem]{Lam2005}. We thus see that rigid fields always have level 1,2, or $\infty$.\\

Recall that a field is called pythagorean if any sum of squares is square. 
Following \cite{ElmanLamSupPyth}, we introduce the following name for formally real rigid fields. 

\begin{KorDef}\label{FormallyRealRigidPythagorean}
	If $F$ is a formally real rigid field, it is pythagorean. A formally real rigid field $F$ is also called \textit{superpythagorean}.
\end{KorDef}
\begin{Beweis}
	If $F$ is formally real and rigid, its Witt ring is isomorphic to $\Z[G]$ for some group $G$ of exponent 2. We thus have $W_tF=\{0\}$ which is equivalent to $F$ being pythagorean by \cite[Chapter VIII., Theorem 4.1 (1)]{Lam2005}.
\end{Beweis}

The above characterization together with Springer's theorem for complete discrete valuation fields motivate us to build the following prototypes of rigid fields in which we can calculate reasonably well and such that these fields realize any possible Witt ring of rigid fields.

\begin{Korollar}\label{AufStandardRigidFieldReduzieren}
	Let $F$ be a rigid field. Then there is a field $K\in\{\mathbb F_3, \R, \C\}$ and an index set $I$ with
	$$WF\cong WK\laurent{t_i}_{i\in I}.$$
\end{Korollar}
\begin{Beweis}
	According to \ref{RigidCharakterisierung} (\ref{EsIstMirEgaaaaaaaal}), we have $WF\cong \Z/n\Z[G]$ for some $n\in\{0,2,4\}$ and some group $G$ of exponent 2. \\
	\begin{minipage}{0.7\textwidth}
			We choose the field $K$ as shown in the adjacent table:
	\end{minipage}\hfill
	\begin{minipage}{0.3\textwidth}
		\begin{tabular}{|c||c|c|c|}
		\hline$n$ & 0 & 2 & 4\\\hline
		$K$ & $\R$ & $\mathbb C$ & $\mathbb F_3$\\\hline
	\end{tabular}
	\end{minipage}\\
	It is well known that we then have $WK\cong \Z/n\Z$. As $G$ is of exponent 2, it can be seen as a vector space over the fields with two elements $\mathbb F_2$ and thus has an $\mathbb F_2$-basis $(g_i)_{i\in I}$ for some index set $I$. We now consider the field $E:=K\laurent{t_i}_{i\in I}$. We then have 
	$$WE\cong\Z/n\Z[G]$$
	as in the proof of \cite[Lemma 1.6]{Ware} (this is essentially a direct limit argument using Springer's Theorem \cite[Chapter VI. Theorem 1.4]{Lam2005}).
\end{Beweis}

The above result further allows us to always work in explicitly given fields if we want to study rigid fields in general. We will fill in the details in the next remark for future reference.

\begin{Bemerkung}\label{ReductionTechnique}
	Because of the Harrison-Cordes Theorem \cite[Chapter XII. Theorem 1.8]{Lam2005}, the study of quadratic forms over rigid fields can be restricted to study quadratic forms over fields of the form $K\laurent{t_i}_{i\in I}$ for a field $K\in\{\mathbb F_3,\R,\C\}$ and some index set $I$, which can be assumed to be well-ordered due to the well-ordering theorem.\\
	If we want to study a concrete form, it is often even possible to only consider the case that $I$ is finite as the direct limit $K\laurent{t_i}_{i\in I}$ can be regarded as the union of the fields $K\laurent{t_{i_1}}\cdots\laurent{t_{i_r}}$ for some $r\in\N_0$ and $i_1,\ldots, i_r\in I$ with $i_1<\ldots<i_r$, see again the proof of \cite[Lemma 1.6]{Ware}. Thus, if a quadratic form $\varphi$ over $E$ is given, we can take any diagonalization of $\varphi$. In this diagonalization, only finitely many Laurent-variables can occur, say these are $t_{j_1},\ldots t_{j_m}$ with $j_1<\ldots<j_m$. Then, $\varphi$ is already defined over $E':=K\laurent{t_{j_1}}\cdots\laurent{t_{j_m}}$ and we can work over this field. For example, the Pfister number of $\varphi$ over $E'$ is bigger than or equal to the Pfister number of $\varphi$ over $E$ as we have $E'\subseteq E$. Thus the task of finding upper bounds for the Pfister numbers over arbitrary rigid fields is reduced to the task of finding upper bounds for the Pfister numbers over fields of the form $K\laurent{t_1}\cdots\laurent{t_n}$ for some $n\in\N$ and $K\in\{\mathbb F_3,\R, \C\}$.
\end{Bemerkung}

The following corollary will be the key idea to determine asymptotic upper bounds for the Pfister numbers. Its proof combines the above theory with the tools that were developed before over fields equipped with a discrete valuation.

\begin{Korollar}\label{KorCDVAufRigid}
	Let $\varphi\in I^nF$ be a quadratic form over some rigid field $F$ that represents 1 and an element ${a\notin \pm D_F(s(F)\times \qf 1)}$, where we interpret $D_F(\infty\times\qf1)$ as $\displaystyle\bigcup_{n\in\N}D_F(n\times\qf1)$. Then there are quadratic forms $\sigma\in I^nF,\tau\in I^{n-1}F$ with $\dim \sigma<\dim\varphi$ and some $t\in F^\ast$ with ${\varphi=\sigma\perp \Pfister t\otimes \tau}$.
\end{Korollar}
\begin{Beweis}
	Using \ref{ReductionTechnique} and \ref{RigidFieldsTauscheBewertungen}, we are reduced to the case where we have $F=K\laurent{t_1}\cdots\laurent{t_n}$ for some $n\in\N$ with $a=t_n$. But then, the assertion readily follows from \ref{UnterteilenInUnimodulareFormen} and \ref{LemmaSplitExactSequence} as both residue class forms for $a=t_n$ are non-hyperbolic by assumption.
\end{Beweis}

We would like to remark that our above result can be applied in particular to rigid fields $F$ with $s(F)=1$. When specialising to the case $n=3$, we get the main results from \cite[Lemma 1.5]{Raczek2013}, the starting point for the calculation of Pfister numbers in the just cited article.

As usual it may be helpful to study the behaviour of a given quadratic form under field extensions. Thus the following result is essential for us.

\begin{Theorem}{\cite[Corollary 2.8]{Ware}}\label{RigidQuadraticExtension}
	Let $F$ be a rigid field and $K/F$ a quadratic field extension. Then $K$ is also a rigid field.
\end{Theorem}

For later reference, we will now discuss the possible diagonalizations of anisotropic binary forms over rigid fields in detail.

\begin{Proposition}\label{DiagonalisierungenRigidBinaer}
	Let $F$ be a rigid field and $\beta=\qf{x,y}$ be an anisotropic binary form over $F$. By abuse of terminology, we say that two diagonalizations of a quadratic form are the same if they only differ by multiplying some entries with a square. We then have one of the following cases:
	\begin{itemize}
		\item $s(F)=1$, $x,y$ represent different square classes and $\qf{x,y}$ and $\qf{y,x}$ are the only diagonalizations of $\beta$;
		\item $s(F)=2$, $x,y$ represent different square classes and $\qf{x,y}$ and $\qf{y,x}$ are the only diagonalizations of $\beta$;
		\item $s(F)=2$, $x,y$ represent the same square classes and $\qf{x,x}$ and $\qf{-x, -x}$ are the only diagonalizations of $\beta$
		\item $s(F)=\infty$, $x,y$ represent different square classes and $\qf{x,y}$ and $\qf{y,x}$ are the only diagonalizations of $\beta$;
		\item $s(F)=\infty$, $x,y$ represent the same square classes and $\qf{x,x}$ is the only diagonalization of $\beta$.
	\end{itemize}
\end{Proposition}
\begin{Beweis}
	We first note that in general, for any $a\in F^\ast$, we cannot have $a$ and $-a$ in the same diagonalization of an anisotropic quadratic form. 
	In the sequel, we use several times the fact that any entry of a diagonalization is represented by the form. Finally, if $x,y$ represent different square classes, we clearly have $D_F(\beta)=\{x,y\}$ because $F$ is rigid.\\
	If we have $s(F)=1$ we have $x=-x$ in $F^\ast/F^{\ast2}$. It is thus clear that $x,y$ have to represent different square classes. As $F$ is rigid we have $D_F(\beta)=\{x,y\}$ and by the above remarks, this case follows.\\
	For $a\in F^\ast$, we have $D_F(\qf{a,a})=\{a,-a\}$ if $s(F)=2$ and $D_F(\qf{a,a})=\{a\}$ if $s(F)=\infty$ by \ref{FormallyRealRigidPythagorean}.
	Thus, if $x,y$ represent different square classes, they both have to occur in any diagonalization of $\beta$. This readily implies that $\qf{x,y}$ and $\qf{y,x}$ are the only diagonalizations of $\beta$ in the respective cases.\\
	So let now $x,y$ represent the same square class.
	If we have $s(F)=2$, it follows by the remarks at the beginning of the proof that $\qf{x,x}$ and $\qf{-x,-x}$ are the only diagonalizations of $\beta$.\\
	Finally, if we have $s(F)=\infty$, \ref{FormallyRealRigidPythagorean} implies that $\qf{x,x}$ is the only diagonalization of $\beta$.
\end{Beweis}

As a corollary, we will now see what makes the theory of quadratic forms over rigid fields much easier than the general case: if one diagonalization of a given form is known, it is easy to determine all the others.

\begin{Korollar}\label{EindeutigDiagRigid}
	Let $\varphi$ be an anisotropic form over a rigid field $F$. If we have $s(F)\in\{1,\infty\}$ the diagonalization of $\varphi$ is unique up to permuting the entries and multiplying them with squares. If we have $s(F)=2$, the diagonalization of $\varphi$ is unique up to permuting the entries, multiplying them with squares and replacing subforms of the form $\qf{x,x}$ for some $x\in F^\ast$ with $\qf{-x,-x}$.
\end{Korollar}
\begin{Beweis}
	It is clear that any of the operations in the statement of the proposition describes isometries of quadratic forms. Further it is well known that two quadratic forms are isometric if and only if they are chain equivalent, see \cite[Chapter I. Chain Equivalence Theorem 5.2]{Lam2005}. The conclusion thus readily follows from \ref{DiagonalisierungenRigidBinaer}.
\end{Beweis}

\begin{Korollar}\label{ValueSetRigidForms}
	Let $\varphi,\psi$ be quadratic forms over a rigid field $F$ such that $\varphi\perp\psi$ is anisotropic. We then have
	$$D_F(\varphi\perp\psi)=\begin{cases}
		D_F(\varphi)\cup D_F(\psi),&\text{ if }s(F)\in\{1,\infty\}\\
		D_F(\varphi)\cup D_F(\psi)\cup\{x\in F^\ast\mid -x\in D_F(\varphi)\cap D_F(\psi)\},&\text{ if }s(F)=2.
	\end{cases}$$
\end{Korollar}
\begin{Beweis}
	It is well known that we have 
	$$D_F(\varphi\perp\psi)=\bigcup_{x\in D_F(\varphi), y\in D_F(\psi)} D_F(\qf{x,y}),$$
	see for example \cite[Chapter I. exercise 20]{Lam2005}. As the elements that are represented by a quadratic form are exactly those that can occur in a diagonalization, the claim now readily follows from \ref{DiagonalisierungenRigidBinaer}.
\end{Beweis}

In the following, we will record some technical results in order to study how hyperbolic planes can occur in the sum of three quadratic forms over rigid fields.

\begin{Lemma}
	Let $F$ be a rigid field and $\varphi_1,\varphi_2,\varphi_3$ be anisotropic quadratic forms over $F$ such that $\varphi_1\perp\varphi_2$ is anisotropic as well. Then $\varphi_1\perp\varphi_2\perp\varphi_3$ is isotropic if and only if one of the following cases occurs:
	\begin{enumerate}
		\item[(1)] at least one of the forms $\varphi_1\perp\varphi_3$ and $\varphi_2\perp\varphi_3$ is isotropic.
		\item[(2)] we have $s(F)=2$ and $D_F(\varphi_1)\cap D_F(\varphi_2)\cap D_F(\varphi_3)\neq\emptyset$.
	\end{enumerate}
\end{Lemma}
\begin{Beweis}
	The form $(\varphi_1\perp\varphi_2)\perp\varphi_3$ is isotropic if and only if there is some ${x\in D_F(\varphi_1\perp\varphi_2)\cap -D(\varphi_3)}$. As we have determined the value set $D_F(\varphi_1\perp\varphi_2)$ in \ref{ValueSetRigidForms}, the claim readily follows by the validity of the following three easy equivalences for some $x$ as above:
	\begin{align*}
		x\in D_F(\varphi_1)&\iff\varphi_1\perp\varphi_3\text{ is isotropic}\\
		x\in D_F(\varphi_2)&\iff\varphi_2\perp\varphi_3\text{ is isotropic}\\
		-x\in D_F(\varphi_1)\cap D_F(\varphi_2)&\iff -x\in D_F(\varphi_1)\cap D_F(\varphi_2)\cap D_F(\varphi_3).
	\end{align*}
\end{Beweis}

\begin{Lemma}\label{LemFindeUnterformen}
	Let $F$ be a rigid field and $\varphi_1, \varphi_2$ be quadratic forms over $F$ such that the orthogonal sum $\varphi_1\perp\varphi_2$ is anisotropic. Further let $\psi\subseteq \varphi_1\perp\varphi_2$ be a subform of $\varphi_1\perp\varphi_2$. Then there are quadratic forms $\psi_1,\psi_2,\psi_3$ over $F$ such that we have $\psi\cong \psi_1\perp\psi_2\perp\psi_3$ and the forms $\psi_1,\psi_2,\psi_3$ fulfil the following:
	\begin{enumerate}[(a)]
		\item\label{LemFindeUnterformenC} $\psi_1\subseteq \varphi_1, \psi_2\subseteq \varphi_2$;
		\item\label{LemFindeUnterformenD} $\big(D_F(\varphi_1)\cup D_F(\varphi_2)\big)\cap D_F(\psi_3)=\emptyset$;
		\item\label{LemFindeUnterformenA} if we have $s(F)\neq 2$, we further have $\psi_3=0$;
		\item\label{LemFindeUnterformenB} for any $x\in F^\ast$, the form $\qf{x,x}$ is not a subform of $\psi_3$.
	\end{enumerate}
\end{Lemma}
\begin{Beweis}
	We prove the assertion by induction on $\dim\psi$, the initial step $\dim\psi=0$ being trivial. We thus assume $\dim\psi>0$ in the following. We will first show that we can decompose $\psi\cong\psi_1\perp\psi_2\perp\psi_3$ such that \ref{LemFindeUnterformenC}, \ref{LemFindeUnterformenD} and \ref{LemFindeUnterformenA} are fulfilled and finally that any such decomposition fulfils \ref{LemFindeUnterformenB} as well.\\
	If we have 
	$$D_F(\psi)\cap \big(D_F(\varphi_1)\cup D_F(\varphi_2)\big)=\emptyset,$$
	we must have $s(F)=2$ by \ref{ValueSetRigidForms} and we can put $\psi_3=\psi$ and $\psi_1=0=\psi_2$. \\
	Otherwise we choose an arbitrary $x\in D_F(\psi)\cap \big(D_F(\varphi_1)\cup D_F(\varphi_2)\big)$ and write ${\psi\cong\qf x\perp \psi'}$ for some suitable form $\psi'$ over $F$. After renumbering we can assume without loss of generality that we have $x\in D_F(\varphi_1)$. In particular there is a form $\varphi_1'$ such that we have $\varphi_1\cong\qf x\perp\varphi_1'$. Using Witt's Cancellation Theorem, we see that $\psi'$ is a subform of $\varphi_1'\perp\varphi_2$.\\
	By induction hypothesis there are quadratic forms $\psi_1'\subseteq \varphi_1', \psi_2'\subseteq\varphi_2$ and $\psi_3'$ with $\big(D_F(\varphi_1')\cup D_F(\varphi_2)\big)\cap D_F(\psi_3')=\emptyset$, such that we have $\psi'\cong\psi_1'\perp\psi_2'\perp\psi_3'$.\\
	We now put
	\begin{align*}
		\psi_1:=\psi_1'\perp\qf x,~~~\psi_2:=\psi_2',~~~\psi_3:=\psi_3'.
	\end{align*}
	Obviously, we have $\psi\cong \psi_1\perp\psi_2\perp\psi_3$ and $\psi_1\subseteq \varphi_1$ und $\psi_2\subseteq\varphi_2$. We will now prove $\big(D_F(\varphi_1)\cup D_F(\varphi_2)\big)\cap D_F(\psi_3)=\emptyset$.\\
	At first, we note that we have
	\begin{align*}
		D_F(\varphi_1)=\begin{cases}
		D_F(\varphi_1')\cup\{x\}, &\text{if } s(F)=1\\
		D_F(\varphi_1')\cup\{x\}, &\text{if } s(F)=2 \text{ and } x\notin D_F(\varphi_1')\\
		D_F(\varphi_1')\cup\{-x\}, &\text{if } s(F)=2 \text{ and } x\in D_F(\varphi_1')\\
		D_F(\varphi_1')\cup\{x\}, &\text{if } s(F)=\infty \text{ and } x\notin D_F(\varphi_1')\\
		D_F(\varphi_1'), &\text{if } s(F)=\infty \text{ and } x\in D_F(\varphi_1').
		\end{cases}
	\end{align*}
	As we have $\big(D_F(\varphi_1')\cup D_F(\varphi_2)\big)\cap D_F(\psi_3')=\emptyset$ by induction hypothesis, the last case is clear. Since $\psi\cong\psi_1\perp\psi_2\perp\psi_3$ with $x\in D_F(\psi_1)$ is anisotropic, we further cannot have $-x\in D_F(\psi_3)$. Thus, the first and the third case are done.\\
	For the remaining two cases, we have to exclude $x\in D_F(\psi_3)$. Assume the contrary. 
	Since we have $\psi_3=\psi_3'$, the induction hypothesis yields $x\notin D_F(\varphi_1')\cup D_F(\varphi_2)$. But $\psi_3=\psi_3'$ is a subform of $\varphi_1'\perp\varphi_2$ so we have $x\in D_F(\varphi_1'\perp\varphi_2)$. As $F$ is rigid, this is only possible if we have $s(F)=2$ and additionally $-x\in D_F(\varphi_1')\cap D_F(\varphi_2)$, see \ref{ValueSetRigidForms}. But this is impossible since then, $\varphi_1=\qf x\perp\varphi_1'$ would be isotropic. Thus \ref{LemFindeUnterformenD} holds.\\
	To prove \ref{LemFindeUnterformenA}, we now assume $s(F)\neq2$. It is then enough to remark that we have $D_F(\varphi_1\perp \varphi_2)=D_F(\varphi_1)\cup D_F(\varphi_2)$ by \ref{ValueSetRigidForms}. Thus the first case in the induction step never occurs and we get $\psi_3=0$ automatically by proceeding as described above.\\
	Finally, for \ref{LemFindeUnterformenB}, we can assume that we have $s(F)=2$ according to \ref{LemFindeUnterformenA}. If we had $\qf{z,z}\subseteq\psi_3$ for some $z\in F^\ast$, we would have $z,-z\in D_F(\psi_3)\subseteq D_F(\psi)\subseteq D_F(\varphi_1\perp\varphi_2)$. 
	As we have 
	$$D_F(\varphi_1\perp\varphi_2)=D_F(\varphi_1)\cup D_F(\varphi_2)\cup\{-x\mid x\in D_F(\varphi_1)\cap D_F(\varphi_2)\}$$
	by \ref{ValueSetRigidForms} this would contradict the fact that we have
	$$\big(D_F(\varphi_1)\cup D_F(\varphi_2)\big)\cap D_F(\psi_3)=\emptyset$$
	and the conclusion follows.
\end{Beweis}

As a strengthening of the above results, we get the following consequence which gives us a precise description of how three quadratic forms over a rigid field have to be related such that their sum has a prescribed Witt index.

\begin{Korollar}\label{WittIndexSumOf3FormsRigid}
	Let $F$ be a rigid field and $\varphi_1,\varphi_2,\varphi_3$ be anisotropic forms over $F$ such that ${\varphi_1\perp\varphi_2}$ is anisotropic as well. Further let $m\in \N$ be an integer. We then have ${i_W(\varphi_1\perp\varphi_2\perp\varphi_3)\geq m}$ if and only if one of the following cases holds:
	\begin{itemize}
		\item we have $s(F)\neq 2$ and there are quadratic forms $\psi_1\subseteq\varphi_1, \psi_2\subseteq\varphi_2$ over $F$ such that 
		\begin{align*}
			\dim(\psi_1\perp\psi_2)&\geq m ~~~\text{ and }~~~-\psi_1\perp-\psi_2\subseteq\varphi_3;
		\end{align*}
		{or}
		\item we have $s(F)=2$ and there are quadratic forms $\psi_1\subseteq\varphi_1,\psi_2\subseteq\varphi_2$ over $F$ and ${x_1,\ldots, x_r\in F^\ast\setminus (D_F(\varphi_1)\cup D_F(\varphi_2))}$ representing pairwise different square classes such that 
		\begin{align*}
			-\psi_1\perp-\psi_2\perp-\qf{x_1,\ldots, x_r}&\subseteq\varphi_3\\
			\qf{x_1,\ldots, x_r}&\subseteq (\varphi_1\perp -\psi_1)_{\text{an}}\perp (\varphi_2\perp-\psi_2)_{\text{an}},\\
			\dim\psi_1+\dim\psi_2+r&\geq m.
		\end{align*}
	\end{itemize}
\end{Korollar}
\begin{Beweis}
	By an easy induction on the integer $m$ using the uniqueness of the Witt decomposition and the anisotropy of $\varphi_1\perp\varphi_2$, we have $i_W(\varphi_1\perp\varphi_2\perp\varphi_3)\geq m$ if and only if there is some quadratic form $\psi$ over $F$ of dimension at least $m$ such that we have $-\psi\subseteq \varphi_3$ and $\psi\subseteq \varphi_1\perp\varphi_2$. \\
	Thus, to show the if part, it is enough to remark that we can choose 
	$$\psi:=\begin{cases}\psi_1\perp\psi_2,&\text{ if }s(F)\neq2\\
	\psi_1\perp\psi_2\perp\qf{x_1,\ldots, x_r},&\text{ if }s(F)=2 
	\end{cases}$$
	as such a form. To show the only if part, let $\psi$ be given as above. We separate the cases $s(F)\neq2$ and $s(F)=2$. If we have $s(F)\neq 2$, \ref{LemFindeUnterformen} yields that we have a decomposition $\psi=\psi_1\perp\psi_2$ and for these $\psi_1,\psi_2$, the requirements are obviously fulfilled.\\
	So let now $s(F)=2$. We apply \ref{LemFindeUnterformen} again and get a decomposition ${\psi=\psi_1\perp\psi_2\perp\psi_3}$, where we can write $\psi_3=\qf{x_1,\ldots, x_r}$ for some ${r\in\N}$ and ${x_1,\ldots, x_r\in F^\ast}$ representing different square classes. As the other properties are readily seen to be satisfied, it remains to show that we have ${\qf{x_1,\ldots, x_r}\subseteq (\varphi_1\perp -\psi_1)_{\text{an}}\perp (\varphi_2\perp-\psi_2)_{\text{an}}}$ As $\psi_i$ is a subform of $\varphi_i$ for $i\in\{1,2\}$ and $\varphi_1\perp\varphi_2$ is anisotropic, the latter form is isometric to ${(\varphi_1\perp\varphi_2\perp-\psi_1\perp-\psi_2)_\text{an}}$. Since we have 
	$$\psi=\psi_1\perp\psi_2\perp\qf{x_1,\ldots, x_r}\subseteq\varphi_1\perp\varphi_2$$
	we get the desired subform relation as an easy consequence of Witt's Cancellation Theorem.
\end{Beweis}

\begin{Lemma}\label{EinWeiteresVielZuTechnischesResultatDasIchNurEinmalBraucheUndDeshalbSoLangGelabeltWerdenKann}
	Let $F$ be a rigid field of level $s(F)=2$ and let $x_1,\ldots, x_r\in F^\ast$ represent pairwise different square classes such that the quadratic form $\qf{x_1,\ldots, x_r}$ is anisotropic. Further, let $\varphi,\psi$ be quadratic forms over $F$ such that $\varphi\perp\psi$ is anisotropic and such that we have $x_i\notin D_F(\varphi)\cup D_F(\psi)$, but $x_i\in D_F(\varphi\perp\psi)$ for all $i\in\{1,\ldots, r\}$. We then have both 
	$$-\qf{x_1,\ldots, x_r}\subseteq\varphi~~~ \text{ and }~~~ -\qf{x_1,\ldots, x_r}\subseteq\psi.$$
\end{Lemma}
\begin{Beweis}
	As we have $\qf{x_1,\ldots, x_r}\subseteq \varphi\perp\psi$ but $x_i\notin D_F(\varphi)\cup D_F(\psi)$ for all $i\in\{1,\ldots, r\}$, \ref{ValueSetRigidForms} implies $-x_i\in D_F(\varphi)\cap D_F(\psi)$. Thus, the induction base is clear by the {Representation Criterion}.
	So let now $r\geq 2$. 
	By the above, we further have representations
	$$\varphi=\qf{-x_1}\perp\varphi'~~~\text{ and }~~~\psi=\qf{-x_1}\perp\psi'.$$
	We thus have 
	$$\varphi\perp\psi\cong\qf{x_1,x_1}\perp\varphi'\perp\psi'$$
	and \ref{ValueSetRigidForms} then implies that we have a disjoint union
	$$D_F(\varphi\perp\psi)=D_F(\qf{x_1,x_1}\perp\varphi'\perp\psi')=\{\pm x_1\}\cup D_F(\varphi'\perp\psi').$$
	Since the form $\qf{x_1,\ldots, x_r}$ is anisotropic and the $x_i$ represent different square classes, we have $x_2,\ldots, x_r\notin\{\pm x_1\}$. We thus have $x_2,\ldots, x_r\in D_F(\varphi'\perp\psi')$.\\
	It is clear that we still have $x_i\notin D_F(\varphi')\cup D_F(\psi')$ for all $i\in\{2,\ldots, r\}$ as these are subforms of $\varphi$ respective $\psi$.
	By induction hypothesis, we have
	$$-\qf{x_2,\ldots,x_r}\subseteq\varphi'~~~\text{ and }~~~-\qf{x_2,\ldots, x_r}\subseteq\psi'$$
	which then implies the assertion.
\end{Beweis}

\section{14-dimensional $I\mathbf{^3}$-forms and 8-dimensional $I\mathbf{^2}$-forms}

From \cite{HoffmannTignol} and \cite{IzhboldinKarpenko}, it is known that there is a deep connection between 14-dimensional $I{^3}$-forms and 8-dimensional $I{^2}$-forms. In this section, we will study both types over rigid fields since the results obtained here will help us to classify 16-dimensional forms in the third power of the fundamental ideal $I^3F$. In this context, it is convenient to introduce the following notation.
\begin{Definition}
	A field $F$ is called a $D(8)$\textit{-field}, if any 8-dimensional form in $I^2F$ whose Clifford invariant has index 4 is Witt equivalent to a sum of 2 forms in $GP_2F$.
	\\The field $F$ is called a $D(14)$\textit{-field} if any 14-dimensional form in $I^3F$ is Witt equivalent to a sum of two forms in $GP_3F$.
\end{Definition}

We will see that rigid fields fulfil both $D(8)$ and $D(14)$. Before proving this, we repeat the classification theorem for 14-dimension $I^3$-forms.

\begin{Proposition}{\cite[Proposition 2.3]{HoffmannTignol} or \cite[Proposition 17.2]{IzhboldinKarpenko}}\label{14DimPfisterZahl}
	Let $\varphi\in I^3F$ be a quadratic form over $F$ with $\dim\varphi=14$. Then $\varphi$ is Witt equivalent to a sum of 3 $GP_3$-forms. Further the following are equivalent:
	\begin{enumerate}[(i)]
		\item there are $\tau_1,\tau_2\in P_3F$ and $s_1, s_2\in F^\ast$ such that $\varphi$ is Witt equivalent to $s_1\tau_1\perp s_2\tau_2$;
		\item there are $\tau_1,\tau_2\in P_3F$ and $s\in F^\ast$ such that $\varphi$ is isometric to $s(\tau_1'\perp-\tau_2')$;
		\item there is some $\sigma\in GP_2F$ with $\sigma\subseteq\varphi$.
	\end{enumerate}
\end{Proposition}

\begin{Proposition}\label{PropRigidD14Field}
	Let $F$ be a rigid field and $\varphi\in I^3F$ be an anisotropic 14-dimensional quadratic form. Then we have $\varphi=\pi_1+\pi_2\in WF$ for some $\pi_1,\pi_2\in GP_3F$, i.e. $F$ is a $D(14)$-field.
\end{Proposition}
\begin{Beweis}
	By \cite[Proposition 2.3]{HoffmannTignol}, we know that $\varphi$ is Witt equivalent to a sum of three forms in $GP_3F$. We choose such a representation $\varphi=\pi_1+\pi_2+\pi_3$ such that we have
	$$\dim(\pi_1\perp\pi_2)_\text{an}\leq\dim(\pi_i\perp\pi_j)_\text{an}$$ for any $i,j\in\{1,2,3\}$ with $i\neq j$. We will distinguish between the possible values that can occur. As we have ${(\pi_1\perp\pi_2)_\text{an}\in I^3F}$, the Gap Theorem implies ${\dim(\pi_1\perp\pi_2)_\text{an}\in\{0,8,12,14,16\}}$.
	\begin{description}
		\item[$\dim=0$:] This contradicts the fact that we have $\dim\varphi=14$.
		\item[$\dim=8$:] In this case, we have $(\pi_1\perp\pi_2)_\text{an}\in GP_3F$ according to the Arason-Pfister Hauptsatz and the claim follows.
		\item[$\dim=12$:] Here Linkage theory implies that $\pi_1$ and $\pi_2$ are both divisible by the same binary Pfister form $\Pfister a$ for some $a\in F^\ast$. In particular $(\pi_1\perp\pi_2)_{F(\sqrt a)}$ is hyperbolic, which then implies 
		$$\dim(\varphi_{F(\sqrt a)})_\text{an}\leq\dim\left((\pi_3)_{F(\sqrt a)}\right)_\text{an}\leq8.$$
		Thus $\varphi$ has a form in $GP_2F$ as a subform. Finally \ref{14DimPfisterZahl} then implies $\varphi$ to be Witt equivalent to a sum of two forms in $GP_3F$.
		\item[$\dim=14$:] According to \ref{14DimPfisterZahl} we can assume, possibly after scaling, that we have $\pi_1,\pi_2\in P_3F$ and $(\pi_1\perp\pi_2)_{\text{an}}=\pi_1'\perp-\pi_2'$, where the prime symbol denotes the pure part of the respective Pfister form as usual.\\
		Further we have $i_W((\pi_1\perp\pi_2)_{\text{an}}\perp\pi_3)=4$. This implies the existence of a quadratic form $\psi$ over $F$ with $\dim\psi=4$, $-\psi\subseteq \pi_3$ and $\psi\subseteq \pi_1'\perp-\pi_2'$. 
		We now decompose $\psi=\psi_1\perp\psi_2\perp\psi_3$ as in \ref{LemFindeUnterformen}. We then have ${\dim\psi_1, \dim\psi_2\leq 1}$ since if we had say $\dim\psi_1\geq2$, we would have
		$$\dim(\pi_1\perp\pi_3)_\text{an}\leq\dim\pi_1+\dim\pi_3-2\dim\psi_1\leq8+8-2\cdot2=12,$$
		contradicting the minimality of $\dim(\pi_1\perp\pi_2)_{\text{an}}$. In particular, we must have $s(F)=2$.\\
		Thus we have $\dim\psi_3\geq2$. 
		According to \ref{LemFindeUnterformen} (\ref{LemFindeUnterformenD}) and (\ref{LemFindeUnterformenB}) there are $${x,y\in D_F(\pi_1'\perp-\pi_2')\setminus(D_F(\pi_1')\cup D_F(-\pi_2'))}$$ that represent different square classes and are represented by $\psi_3$. Now \ref{EinWeiteresVielZuTechnischesResultatDasIchNurEinmalBraucheUndDeshalbSoLangGelabeltWerdenKann} implies 
		$$-\qf{x,y}\subseteq\pi_1'~~~\text{ and }~~~\qf{x,y}\subseteq\pi_2'.$$
		This implies that both $\pi_1$ and $\pi_2$ become isotropic (hence hyperbolic) over $F(\sqrt{-xy})$. Since this is equivalent to $\pi_1,\pi_2$ having a common slot and as we have ${\dim(\pi_1\perp\pi_2)_\text{an}=14}$, this contradict Linkage theory.
		\item[$\dim=16$:] Just as above in the case $\dim(\pi_1\perp\pi_2)_\text{an}=14$, we can deduce that the Pfister forms that $\pi_1$ respectively $\pi_2$ are similar to have a common slot. Thus, as in the case $\dim(\pi_1\perp\pi_2)_\text{an} =12$, we see that $\varphi$ contains a subform in $GP_2F$ and is thus Witt equivalent to a sum of two forms in $GP_3F$ according to \ref{14DimPfisterZahl} again.
	\end{description}
\end{Beweis}

Because of the strong connection of the two types of forms studied here, we can easily deduce the following as announced before:

\begin{Korollar}\label{RigidD8Field}
	Rigid fields are $D(8)$-fields.
\end{Korollar}
\begin{Beweis}
	Since $F$ is rigid, so is $F\laurent t$ according to \cite[Examples 1.11 (iv)]{Ware}. As we have shown in \ref{PropRigidD14Field}, $F\laurent t$ is a $D(14)$-field. By \cite[Theorem 4.1]{HoffmannTignol}, this implies $F$ to be a $D(8)$-field.
\end{Beweis}

It would be interesting to prove $D(8)$ directly, such that we can get $D(14)$ by {\cite[Theorem 4.4]{HoffmannTignol}}.

\section{16-dimensional $I\mathbf{^3}$-forms}

We are able to classify those 16-dimensional forms in $I^3F$ for rigid fields that are Witt equivalent to a sum of at most three forms in $GP_3F$. Its proof uses the same techniques as the proof of \ref{PropRigidD14Field}. At the end of the section, we will see that any 16-dimensional form in $I^3F$ satisfies the following equivalent conditions. 

\begin{Proposition}\label{RigidCharacterization16Dim}
	Let $F$ be a rigid field and $\varphi\in I^3F$ be an anisotropic quadratic form with $\dim\varphi=16$. Then the following are equivalent:
	\begin{enumerate}[(i)]
		\item $\varphi$ is isometric to a sum of 4 forms in $GP_2F$;
		\item $\varphi$ contains a subform in $GP_2F$;
		\item $\varphi$ is Witt equivalent to a sum of at most at most 3 forms in $GP_3F$.
	\end{enumerate}
\end{Proposition}
\begin{Beweis}
	The implication (i) $\Rightarrow$ (ii) is trivial. For the implication (ii) $\Rightarrow$ (iii), we write $\varphi=\sigma\perp \qf w\perp\psi$ for some $\sigma\in GP_2F$, some $w\in F^\ast$ and a suitable quadratic form $\psi$ of dimension 11 over $F$. We can find $x,y,z\in F^\ast$ such that we have $\sigma\perp\qf{w,x,y,z}\in GP_3F$. In $WF$ we thus have
	\begin{align*}
		\varphi=\sigma\perp \qf w\perp\psi=\sigma\perp \qf{w,x,y,z}\perp\psi\perp\qf{-x,-y,-z}.
	\end{align*}
	We further have $\dim(\psi\perp\qf{-x,-y,-z})=14$ so that this form is Witt equivalent to a sum of at most two $GP_3F$-forms by \ref{PropRigidD14Field}, so this implication is done.
	
	For the implication (iii) $\Rightarrow$ (i), let now $\varphi=\pi_1+\pi_2+\pi_3\in WF$ with $\pi_1,\pi_2,\pi_3\in GP_3F$. We further assume that $\dim(\pi_1\perp\pi_2)_\text{an}$ is minimal under all such representations.
	
Similar to the proof of \ref{PropRigidD14Field}, we are readily reduced to the cases in which we have $\dim(\pi_1\perp\pi_2)_\text{an}\in\{8, 12, 14, 16\}$.\\
	$\dim=8$:
	If we have $\dim(\pi_1\perp\pi_2)_\text{an}=8$ then $(\pi_1+\pi_2)_{\text{an}}$ is isometric to some $\pi\in GP_3F$ according to the Arason-Pfister Hauptsatz. Thus, we have $\varphi\cong\pi\perp\pi_3$ and the conclusion follows.
	\\
	$\dim=12$:
	If we have $\dim(\pi_1\perp\pi_2)_\text{an}=12$, then $(\pi_1+\pi_2)_{\text{an}}$ is divisible by a binary form $\Pfister a$ due to \cite[Satz 14, Zusatz]{Pfister}. Thus, we have $i_W(\varphi_{F(\sqrt a)})\geq 4$ and we can write $\varphi\cong\Pfister{a}\otimes \sigma\perp\psi$ with some $4$-dimensional form $\sigma$ and some $8$-dimensional form $\psi$. According to \cite[Example 9.12]{Knebusch2} $\Pfister{a}\otimes \sigma$ is an 8-dimensional form in $I^2F$, whose Clifford invariant has index at most $2$. In $WF$ we therefore have
	$$\psi=\varphi-\Pfister a\otimes\sigma\in I^2F$$
	which then implies
	$$c(\psi)=c(\varphi)c(\Pfister a\otimes\sigma)=c(\Pfister a\otimes\sigma).$$
	Using \cite[Example 9.12]{Knebusch2} again, we see that $\psi$ is divisible by a binary form as well. As 8-dimensional forms that are divisible by a binary form are isometric to a sum of two forms in $GP_2F$, we are done in this case.
	\\
	$\dim=14$:
	So let now $\dim(\pi_1\perp\pi_2)_{\text{an}}=14$. According to \ref{14DimPfisterZahl} we can assume we have $(\pi_1\perp\pi_2)_{\text{an}}\cong\pi_1'\perp-\pi_2'$, possibly after a scaling. We further have ${i_W(\pi_1'\perp-\pi_2'\perp\pi_3)=3}$, such that there is some 3-dimensional form $\psi\subseteq \pi_1'\perp-\pi_2'$ with $-\psi\subseteq \pi_3$. 
	
	We decompose $\psi=\psi_1\perp\psi_2\perp\psi_3$ as in \ref{LemFindeUnterformen}. Because of the minimality of $\dim(\pi_1\perp\pi_2)_{\text{an}}$, we have $\dim\psi_1\leq 1$ and $\dim\psi_2\leq1$.
	
	As in the case of dimension 14 in the proof of \ref{PropRigidD14Field}, we can see that $\dim\psi_3\geq2$ would contradict Linkage theory. As the dimensions of $\psi_1,\psi_2$ and $\psi_3$ have to sum up to 3, we get
	\[
		\dim\psi_1=\dim\psi_2=\dim\psi_3=1.
	\]
	Thus $\varphi$ contains a 10-dimensional subform that is the orthogonal sum of a 5-dimensional subform of $\pi_1$ and a 5-dimensional subform of $\pi_2$. Both of these forms are Pfister neighbors that contain a subform in $GP_2F$ according to \cite[Chapter X. Proposition 4.19]{Lam2005}.
	Thus $\varphi$ has a decomposition $\varphi=\sigma\perp\tau$, where $\sigma$ is isometric to a sum of 2 forms in $GP_2F$. We thus have $\sigma\in I^2F$ and the Clifford invariant of $\sigma$ has index at most $4$. 
	As in the case $\dim=12$, these properties also hold for $\tau$. Applying \ref{RigidD8Field} now gives us that $\tau$ is also isometric to a sum of two forms in $GP_2F$ which finishes this case.
		\\
		$\dim=16$:
		Here, we are reasoning just as in the latter case and use the same terminology for all upcoming forms etc. We have $\dim\psi=4$. Because of the minimality of $\dim(\pi_1\perp\pi_2)_{\text{an}}$, we even have $\psi_1=0=\psi_2$.
		As in the case $\dim=14$ above (i.e. as in the proof of \ref{PropRigidD14Field}), we see that the Pfister forms that are similar to $\pi_1$ respectively $\pi_2$ have a common slot, so that $\pi_1\perp\pi_2$ is divisible by a binary form $\Pfister a$. Now the conclusion follows as in the case $\dim=12$.
\end{Beweis}

Our next goal is to study 16-dimensional form in $I^3F$ in more detail in order to prove that each such form satisfies the equivalent conditions of \ref{RigidCharacterization16Dim}. To do so, we need the next technical lemma.

\begin{Lemma}\label{RigidLemma14DimTechnisch}
	Let $F$ be a rigid field and $\varphi_1,\varphi_2$ be two anisotropic quadratic forms over $F$, such that $\varphi_1\perp\varphi_2$ is an anisotropic form in $I^3F$ of dimension 14. Then, for any $t\in F^\ast$, the form $\varphi_1\perp t\varphi_2$ contains a subform in $GP_2F$.
\end{Lemma}
\begin{Beweis}
	We show that one of the forms $\varphi_1$ and $\varphi_2$ already contains a subform in $GP_2F$ or that there is some binary form that is similar to both a subform of $\varphi_1$ and a subform of $\varphi_2$. This obviously implies the assertion.
	
	Since $F$ is a rigid field, $F$ is a $D(14)$-field by \ref{PropRigidD14Field}. Therefore, after a possible scaling, we may assume that we have $\pi_1,\pi_2\in P_3F$ with 
	$$\varphi_1\perp\varphi_2\cong \pi_1'\perp-\pi_2'.$$ 
	We remark that $\pi_1,\pi_2$ cannot have a common slot.
	
	As $\pi_1,\pi_2$ are 3-fold Pfister forms, we can choose $a, a'\in F^\ast$ and 3-dimensional forms $\sigma,\sigma'$ over $F$ such that we have 
	\begin{align*}
		\psi:=\Pfister a\otimes \sigma\subseteq \pi_1', ~~ \psi':=\Pfister{a'}\otimes \sigma'\subseteq -\pi_2',
	\end{align*}
	see \cite[Theorem 4.1]{HoffmannSplitting}.
	In particular $\psi\perp\psi'$ is also a subform of $\varphi_1\perp\varphi_2$. We now decompose $\psi\cong\psi_1\perp\psi_2\perp\psi_3$ and $\psi'\cong\psi_1'\perp\psi_2'\perp\psi_3'$ according to \ref{LemFindeUnterformen}. 
	We will now proof the assertion while distinguishing the possible dimensions of these subforms:\\
		\textbf{Case 1:} $\dim\psi_3=0$ or $\dim\psi_3'=0$:~\\
		According to the symmetry of the statement, it is enough to consider the case $\dim\psi_3=0$. Further we can assume $\dim\psi_1\geq\dim\psi_2$, possibly after renumbering the $\varphi_i$. As we clearly have $\dim\psi_1+\dim\psi_2=6$ the latter implies $\dim\psi_1\geq 3$. 
		 
		 If we have $\dim\psi_1\geq 5$, it follows, readily that $\psi_1$ already contains a four dimensional subform that is divisible by $\Pfister a$, i.e. a form in $GP_2F$.
		 
		 If we have $\dim\psi_1=4$ we can use the same arguments as above to get that $\psi_1$ becomes isotropic over $F(\sqrt a)$ which then implies that $\psi_1$ is similar to $\Pfister a\perp\sigma$ with some quadratic form $\sigma$ of dimension 2. We then have that $\sigma\perp\psi_2$ is divisible by $\Pfister a$, i.e. a form in $GP_2F$. Using \cite[Chapter X. Corollary 5.4]{Lam2005} one readily sees that this is only possible if $\sigma$ and $\psi_2$ are similar which concludes this case. 
		 
		 If $\dim\psi_1=3$ and $\psi_1$ becomes isotropic over $F(\sqrt a)$, then so does $\psi_2$ as $\psi$ becomes hyperbolic over $F(\sqrt a)$. Thus, both $\psi_1$ and $\psi_2$ contain a subform similar to $\Pfister a$ and this case is done.
		 \\
		 Otherwise $\psi_1$ and $\psi_2$ are quadratic forms of dimension 3 that stay anisotropic over $K:=F(\sqrt a)$ but fulfil $(\psi_1)_K\cong-(\psi_2)_K$. 
		 By \ref{RigidQuadraticExtension} $K$ is a rigid field, too. Using \ref{EindeutigDiagRigid} we see that the diagonalization of $(\psi_1)_K$ is either unique up to multiplying its entries with squares and permuting the entries or we have $s(K)=2$ (and thus also $s(F)=2$ as can readily seen using \cite[Theorem 2.7]{Ware}) and $(\psi_1)_K=\qf{x,x,y}$ for some $x,y\in F^\ast$. 
		 \\
		 In the first case, we write $(\psi_1)_K=\qf{x,y,z}$ for suitable $x,y,z\in F^\ast$ representing pairwise different square classes in $K$. Using \cite[Chapter VII. Theorem 3.8]{Lam2005}, we see that we have
		 $$	\psi_1=\qf{a^{i_1}x, a^{j_1}y, a^{k_1}z}\text{ and }	\psi_2=-\qf{a^{i_2}x, a^{j_2}y, a^{k_2}z}$$
		 for some $i_1,i_2, j_1, j_2, k_1, k_2\in\{0, 1\}$. After renaming $x,y,z$, the pigeon hole principle implies that we have either $i_1=i_2$ and $j_1=j_2$ or $i_1\neq i_2$ and $j_1\neq j_2$. In both cases $\qf{a^{i_1}x, a^{j_1}y}$ and $-\qf{a^{i_2}x, a^{j_2}y}$ are similar so that this case is done.\\
		 In the second case we argue the same way. We get that $\psi_1$ is isometric to one of the following forms on the left for some $i\in \{0,1\}$ and $\psi_2$ is isometric to one of the forms on the right for some $j\in\{0,1\}$: \\
		 \begin{minipage}{0.5\textwidth}
			\begin{align*}
		 		\qf{x,x,a^iy}\cong \qf{-x,-x, a^iy}\\
				\qf{ax, ax, a^iy}\cong\qf{-ax, -ax, a^iy}\\
				\qf{-x, -ax, a^iy}\\
				\qf{x, ax, a^{i}y}
			\end{align*}
		\end{minipage}
		\begin{minipage}{0.5\textwidth}
			\begin{align*}
		 	\qf{x,x,-a^jy}\cong \qf{-x,-x, -a^jy}\\
			\qf{ax, ax, -a^jy}\cong\qf{-ax, -ax, -a^jy}\\
			\qf{-x, -ax, -a^jy}\\
			\qf{x, ax, -a^{j}y}
		\end{align*}
		\end{minipage}
	
		Thus a binary form that is similar to both a subform of $\psi_1$ and a subform of $\psi_2$ can be found in the upcoming table in which all cases with $\psi_1\not\cong-\psi_2$ (that case being clear) are considered.
		
		\begin{tabular}{|c||c|c|c|c|}
			\hline
			& $\qf{x,x,-ay}$ & $\qf{ax, ax, -ay}$ & $\qf{-x, -ax, -ay}$ & $\qf{x, ax, -ay}$ \\\hline\hline
			$\qf{x,x,y}$ & $\qf{x,x}$ & $\qf{x,x}$ & $\qf{x,y}$ & $\qf{-x,y}$\\ \hline
			$\qf{ax, ax, y}$ & $\qf{x,x}$ & $\qf{x,x}$ & $\qf{ax,y}$ & $\qf{-ax,y}$\\ \hline
			$\qf{-x, -ax, y}$ & $\qf{-ax,y}$ & $\qf{-x,y}$ & $\qf{-x,-ax}$ & $\qf{-x,-ax}$\\ \hline
			$\qf{x, ax, y}$ & $\qf{ax,y}$ & $\qf{x,y}$ & $\qf{x,ax}$ & $\qf{x,ax}$\\ \hline
		\end{tabular}
		\\
		
		\textbf{Case 2:} $\dim\psi_3\geq2$ or $\dim\psi_3'\geq2$:~\\
		 	It is again enough to consider the case $\dim\psi_3\geq 2$. Because of \ref{LemFindeUnterformen} (\ref{LemFindeUnterformenB}) there are $x,y\in F^\ast$ representing different square classes with $\psi_3=\qf{x,y,\ldots}$. Because of \ref{LemFindeUnterformen} \ref{LemFindeUnterformenD} we have $x,y\in D_F(\psi_3)\subseteq D_F(\varphi_1\perp\varphi_2)$ but $x,y\notin D_F(\varphi_1)\cup D_F(\varphi_2)$. Now, \ref{EinWeiteresVielZuTechnischesResultatDasIchNurEinmalBraucheUndDeshalbSoLangGelabeltWerdenKann} implies both $\varphi_1=\qf{-x,-y,\ldots}$ and $\varphi_2=\qf{-x,-y,\ldots}$. According to the statement at the beginning of the proof, this case is done.\\
			
		\textbf{Case 3:} $\dim\psi_3=1=\dim\psi_3'$:\\
			If we have $\psi_3=\qf x\not\cong\qf y=\psi_3'$ for some $x,y\in F^\ast$, we can argue as in the last case using \ref{EinWeiteresVielZuTechnischesResultatDasIchNurEinmalBraucheUndDeshalbSoLangGelabeltWerdenKann} to get $\varphi_1=\qf{-x,-y,\ldots}, \varphi_2=\qf{-x,-y,\ldots}$ and we are done.\\
			Otherwise we have $\psi_3=\qf x=\psi_3'$, so we can write ${\varphi_1=\nu_1\perp\qf{-x}}$ and ${\varphi_2=\nu_2\perp\qf{-x}}$. We further choose orthogonal complements of $\qf x$ in $\pi_1'$ respectively $-\pi_2'$. As in the beginning of the proof, we can write them as a product of a Pfister form and a ternary form, i.e. we have
			$$\pi_1'=\Pfister b\otimes \tau\perp\qf{x} \text{   and   } -\pi_2'=\Pfister{b'}\otimes\tau'\perp\qf{x}$$
			for some ternary forms $\tau,\tau'$ and $b,b'\in F^\ast$. We have a chain of isometries
			\begin{align*}
				\nu_1\perp\qf x\perp\nu_2\perp\qf x&\cong \nu_1\perp\qf{-x}\perp\nu_2\perp\qf{-x}\\
				&\cong \varphi_1\perp\varphi_2\\
				&\cong\pi_1'\perp-\pi_2'\\
				&\cong \Pfister b\otimes \tau\perp\qf{x}\perp \Pfister{b'}\otimes\tau'\perp\qf{x}.
			\end{align*}
			Witt's cancellation law now implies $\Pfister b\otimes \tau\perp\Pfister{b'}\otimes\tau'\cong \nu_1\perp\nu_2$.
		
			We now apply the above argument for $\Pfister b\otimes \tau$ and $\Pfister{b'}\otimes\tau'$ as subforms of $\nu_1\perp\nu_2$. Note that all arguments used above stay valid as we did not use any specific information on $\varphi_1,\varphi_2$ but only of the chosen subforms $\psi,\psi'$.
			If we are in case 1 or 2 for $b,b',\tau,\tau', \nu_1,\nu_2$ we are done as we have already seen. If we are again in case 3 for $b,b',\tau,\tau', \nu_1,\nu_2$, we get the existence of some $y\in F^\ast$ represented by both $\pi_1'$ and $-\pi_2'$. This would imply $\pi_1$ and $\pi_2$ to have $-xy$ as a common slot similar as in the case $\dim =14$ in \ref{PropRigidD14Field}, which we excluded at the beginning of the proof. Thus we are done.
\end{Beweis}

\begin{Theorem}\label{RigidGP316}
	Let $F$ be a rigid field and $\varphi\in I^3F$ be an anisotropic quadratic form over $F$ of dimension 16. Then $\varphi$ is Witt equivalent to a sum of at most three forms in $GP_3F$.
\end{Theorem}
\begin{Beweis}
	We will show that $\varphi$ contains a subform in $GP_2F$ so that the conclusion then follows by \ref{RigidCharacterization16Dim}. After scaling, we can assume $1\in D_F(\varphi)$. If $\varphi$ is isometric to $16\times \qf{1}$ (which is only possible if $F$ is superpythagorean), the assertion is clear. Otherwise there is some $n\in\N$ such that we can assume $\varphi$ to be defined over the field $K\laurent{t_1}\cdots\laurent{t_n}$ and that $\varphi$ has a decomposition into residue class forms ${\varphi\cong \varphi_1\perp t_n\varphi_2}$ such that both residue class forms have positive dimension. As mentioned in \ref{RigidFieldsTauscheBewertungen} we can replace the uniformizer $t_n$ with $at_n$ for any $a\in K\laurent{t_1}\cdots\laurent{t_{n-1}}^\ast$. By doing so, we also get $a\varphi_2$ as the second residue class form instead of $\varphi_2$.
	 We may thus assume $D_F(\varphi_1)\cap D_F(\varphi_2)\neq\emptyset$, i.e. $\sigma:=(\varphi_1\perp-\varphi_2)_\text{an}$ has dimension at most 14. If we have $\dim\sigma\leq 12$, there is some binary form $\beta$ that is a subform of both $\varphi_1$ and $\varphi_2$, so that $\beta\otimes\qf{1,t_n}\in GP_2F$ is a subform of $\varphi$.\\
	If we have $\dim\sigma=14$, there is some $x\in F^\ast$ and quadratic forms $\psi_1, \psi_2$ such that we have
	$$\varphi_1\cong\qf{x}\perp\psi_1\text{   and   }\varphi_2\cong\qf{x}\perp\psi_2.$$
	As in the proof of \ref{UnterteilenInUnimodulareFormen}, we have $\sigma \cong \psi_1\perp-\psi_2\in I^3F$ (in fact, our $\sigma$ here has exactly the same role as the $\sigma$ in the above mentioned result). As we have $\dim\sigma=14$, it contains a subform lying in $GP_2F$ according to \ref{PropRigidD14Field}. By \ref{RigidLemma14DimTechnisch} the form $\psi_1\perp t_n\psi_2$ also contains a $GP_2$-subform, which then trivially implies
	$$\varphi\cong \psi_1\perp\qf x\perp t_n(\psi_2\perp\qf x)$$
	to have a subform in $GP_2F$, which concludes the proof.
\end{Beweis}

\begin{Beispiel}
	The bound in \ref{RigidGP316} is sharp as the following example shows. Let $K\in\{\mathbb F_3, \R, \C\}$ and $F=K\laurent a\laurent b\laurent c\laurent d\laurent e\laurent f$. We first construct an 8-dimensional form in $I^2F$ that is not Witt equivalent to a sum of $2$ forms in $GP_2F$. To do so, we can consider
	$$\psi:=\qf{1, a, b, c, d, e, f, abcdef}\in I^2F,$$
	which is the generic 8-dimensional form in $I^2F$ and fulfills ${GP}_2(\psi)=3$ by \ref{Maximal2PfisterNumber}. Then, $\varphi:=\psi\otimes\Pfister{t}\in I^3F\laurent t$ fulfils $GP_3(\varphi)=3$ by \ref{PropPfisterZahlUnimodularMalUniformizer}.
\end{Beispiel}

Another common way to measure the complexity of a quadratic form is to study its splitting behaviour over multiquadratic field extensions. There are 16-dimensional $I^3$-forms over non-rigid fields that do not split over multiquadratic extensions of degree $\leq 8$, see \cite[Theorem 2.1]{Karpenko2017}. For rigid fields, the situation is much less involved.

\begin{Proposition}\label{MultiQuadraticSplitting16Rigid}
	Let $\varphi$ be a 16-dimensional form in $I^3F$ with $F$ rigid. Then $\varphi$ splits over some biquadratic extension of $F$, i.e. there are $a,b\in F^\ast$ such that $\varphi_{F(\sqrt a,\sqrt b)}$ is hyperbolic.
\end{Proposition}
\begin{Beweis}
	According to \ref{RigidGP316} and \ref{RigidCharacterization16Dim} we can write $\varphi=\psi\perp\sigma$ where we have $\sigma\in GP_2F$. We choose $a\in F$ such that $\sigma_{F(\sqrt a)}$ is isotropic hence hyperbolic. If $\psi_{F(\sqrt a)}$ is isotropic then it is hyperbolic or Witt equivalent to a form in $GP_3F(\sqrt a)$ that is defined over $F$ as quadratic extensions are excellent, see \cite[Chapter XII. Proposition 4.4]{Lam2005}. In both cases the assertion is clear.\\
	Otherwise $\psi_{F(\sqrt a)}$ is an anisotropic, 12-dimensional form in $I^3F(\sqrt a)$ and hence divisible by a binary Pfisterform $\Pfister b$ for some $b\in K^\ast$. By \cite[Theorem 1.9]{Ware}, the square class of $b$ in $F(\sqrt a)$ has a representative of the form $z$ or $z\sqrt a$ for some $z\in F^\ast$. We are done if we can exclude the latter case. As $F(\sqrt a)$ is also a rigid field by \ref{RigidQuadraticExtension}, 
	we know how two diagonalizations of the same form can differ by \ref{EindeutigDiagRigid}. As $\psi$ is defined over $F$, we can thus deduce that we must have $b\in F^\ast$.
\end{Beweis}

\begin{Beispiel}
	\ref{MultiQuadraticSplitting16Rigid} is sharp in the sense that in general, forms over dimension 16 in $I^3F$ over a rigid field $F$ will not split over a quadratic extension. As an example, we can consider the 16-dimensional form 
	${\Pfister{a,b,c}\perp\Pfister{d,e,f}}$ over the field ${F:=K\laurent a\laurent b\laurent c\laurent d\laurent e\laurent f}$ where we can choose $K\in\{\R, \C, \mathbb F_3\}$. 
\end{Beispiel}

We can show, that the characterization in \ref{RigidCharacterization16Dim} does not generalize to arbitrary fields. To do so, we need the following result.

\begin{Proposition}\label{PropPfisterZahl16BeiUnterform}
	Let $\varphi\in I^3F$ be an anisotropic quadratic form with $\dim\varphi=16$. We further presume the existence of some $\sigma,\tau\in GP_2F$ with $\sigma\perp\tau\subseteq \varphi$. Then $\varphi$ is Witt equivalent to a sum of at most three elements in $GP_3F$. 
\end{Proposition}
\begin{Beweis}
	By our assumption, we have $\varphi\cong \sigma\perp\tau\perp\qf{w}\perp\psi$ for some $w\in F^\ast$ and a 7-dimensional quadratic form $\psi$ over $F$. We choose $x,y,z\in F^\ast$ such that $\qf{w,x,y,z}$ is similar to $\sigma$. This implies in particular $\sigma\perp\qf{w,x,y,z}\in GP_3F$. In $WF$ we thus have
	\begin{align*}
		\varphi=\sigma+\tau+\qf{w}+\psi=\left(\sigma\perp\qf{w,x,y,z}\right)+\left(\tau\perp\qf{-x,-y,-z}\perp\psi\right).
	\end{align*}
	Since we have $\varphi, \sigma\perp\qf{w,x,y,z}\in I^3F$, we also have $\tau\perp\qf{-x,-y,-z}\perp\psi\in I^3F$. Further we have $\dim(\tau\perp\qf{-x,-y,-z}\perp\psi)= 14$ and this form contains $\tau\in GP_2F$ as a subform. Thus $\tau\perp\qf{-x,-y,-z}\perp\psi$ is Witt equivalent to a sum of at most two $GP_3F$-forms by \ref{14DimPfisterZahl} and the conclusion follows.
\end{Beweis}

\begin{Beispiel}\label{BeispielAufDasIchVorherSchonMalVerweise}
	In order to show that the characterization in \ref{RigidCharacterization16Dim} does not hold over non-rigid fields, we will construct a 16-dimensional form in $I^3F$ for a suitable field $F$ that has Pfister number 3 but is not isometric to a sum of four forms in $GP_2F$. Over the field $F:=\Q(x)\laurent{t_1}\cdots\laurent{t_4}$ we consider the forms 
	\begin{align*}
		\psi_1&:=\qf{x,-(x+4)}\perp-t_1\qf{1, -(x+4)},\\
		\psi_2&:=\qf{x, -(x+1)}\perp -2t_1\qf{1, -(x+1)}\\
		\rho_1&:=\qf{1, -x, -t_1t_2(x+4), t_1t_2x(x+4)}=\Pfister{x, t_1t_2(x+4)},\\
		\rho_2&:=\qf{1, -x(x+1)(x+4), 2t_1x(x+2), -2t_1(x+1)(x+2)(x+4)}\\
		&=\Pfister{x(x+1)(x+4), -2t_1x(x+2)}.
	\end{align*}
	and finally build the form $\varphi:=\psi_1\perp-t_2\psi_2\perp t_4(\rho_1\perp t_3\rho_2)$. In the sequel we will use a lot of facts shown in \cite[Example 6.3]{HoffmannTignol}. At first we know that $\varphi_1:=\psi_1\perp-t_2\psi_2$ is anisotropic, lies in $I^2F$ and does not contain a subform in $GP_2F$. Further, $\varphi_2:=\rho_1\perp t_3\rho_2$ is also an anisotropic form in $I^2F$ that has the same Clifford invariant. We thus have $\varphi\in I^3F$ with $\dim\varphi=16$ and $\varphi$ is anisotropic. By \ref{PropPfisterZahl16BeiUnterform} we further know that $\varphi$ has 3-Pfister number at most 3. By showing that $\varphi$ is not isometric to a sum of four forms in $GP_2F$, it will further be clear that we even have an equality $GP_3(\varphi)=3$.\\
	Similarly as in \ref{LemFindeUnterformen}, we can show by an induction argument that any form $\psi\subseteq\varphi$ has a decomposition $\psi\cong\psi_1\perp t_4\psi_2$ with $\psi_1\subseteq\varphi_1, \psi_2\subseteq\varphi_2$. Thus, if $\varphi$ would be isometric to an orthogonal sum of four $GP_2$-forms, there has to be a $\sigma\in GP_2F$ that can be decomposed into $\sigma_1\perp t_4\sigma_2$ with $\sigma_1\subseteq\varphi_1,\sigma_2\subseteq\varphi_2$ and $\sigma_1\neq0\neq\sigma_2$ (as $\varphi_1$ does not contain any subform in $GP_2F$ itself). By \ref{LemmaSplitExactSequence} we have $\sigma_1,\sigma_2\in IF$ which then implies $\dim\sigma_1=\dim\sigma_2=2$. As we have $\det\sigma=1$, we have $\det\sigma_1=\det\sigma_2$. Analysing the decomposition of $\varphi$, we see this can only happen if we have $\sigma_2\subseteq\rho_1$ or $\sigma_2\subseteq t_3\rho_2$. We will assume $\sigma_2\subseteq \rho_1$, the other case is similar.\\
	We now choose an $a\in F^\ast$ such that $\sigma_2$ becomes isotropic (hence hyperbolic) over $F(\sqrt a)$. In fact, by the choice of $\sigma_2$, we can even choose $a\in \Q(x)\laurent{t_1}\laurent{t_2}^\ast$. Then, as Pfister forms are either anisotropic or hyperbolic and $\sigma_1$ is similar to $\sigma_2$, both $\sigma_1$ and $\rho_1$ become hyperbolic over $F(\sqrt a)$. This implies $i_W(\varphi_{F(\sqrt a)})\geq 3$ and thus, by the Gap Theorem, even $i_W(\varphi_{F(\sqrt a)})\geq4$. By the choice of $a$, the $t_4$-adic valuation has an extension to $F(\sqrt a)$ and $t_4$ still is a uniformizer. As the Witt index of a form over a complete discrete valuation field is the sum of the Witt indices of its residue class forms (even though they might not be unique up to isometry), the inequality $i_W(\varphi_{F(\sqrt a)})\geq4$ can only be fulfilled if 
	\begin{enumerate}[(a)]
		\item\label{LabelBspA} $(\rho_2)_{F(\sqrt a)}$ is isotropic or 
		\item\label{LabelBspB} we have $i_W((\varphi_1)_{F(\sqrt a)})\geq 2$.
	\end{enumerate}
	In case (\ref{LabelBspA}), the Pfister forms that $\rho_1$ respectively $\rho_2$ are similar to would have a common slot, but this was excluded in \cite[Example 6.3]{HoffmannTignol}.\\
	But case (\ref{LabelBspB}) would imply the existence of a subform of $\varphi_1$ lying in $GP_2F$, a contradiction. Thus the proof is complete.
\end{Beispiel}

\section{Asymptotic Pfister Numbers}

In this section, we will study the growth of Pfister numbers for forms of increasing dimension. As a fixed field can be too small to have anisotropic forms of all dimensions, which is a necessary assumption to talk about meaningful lower bounds, we will allow rigid field extensions while finding lower bounds as can be seen in the upcoming Proposition.

\begin{Proposition}
	Let $F$ be a rigid field. Then, there is some field extension $E/F$ such that $E$ is a rigid field and for any integer $d\geq8$, we have
	\begin{align}\label{MedicalDetectives}
		GP_3(E,d)\geq \left\lfloor \frac d4\right\rfloor-1.
	\end{align}
\end{Proposition}
\begin{Beweis}
	As the term on the right sight of \eqref{MedicalDetectives} increases monotonously when $d$ grows, we may assume that $d$ is even.
	According to \ref{AufStandardRigidFieldReduzieren} and passing to a field extension, we may further assume ${F=K\laurent{t_i}_{i\in I}}$ for some algebraically closed field $K$ and some infinite index set $I$. To simplify notation, we assume $\N\subseteq I$. 
	We define the integer $n$ to be
	$$n:=2\cdot\left\lfloor\frac d4\right\rfloor-2 = \begin{cases}
		\frac d2-2,&\text{ if }d\equiv0\mod4\\
		\frac{d}2-3,&\text{ if }d\equiv2\mod4
	\end{cases}.$$
	Note that $n$ is even in both cases.
	By \ref{Maximal2PfisterNumber}, using \ref{UeberLaurentBleibtGleich} and induction (recall the definition of $K\laurent{t_i}_{i\in I}$ as a direct limit, see \ref{AufStandardRigidFieldReduzieren} again), for $\psi:=\qf{1, t_1,\ldots, t_{n},(-1)^{\frac{n+2}2} t_1\cdot\ldots\cdot t_{n}}\in I^2F$, we have
	$$GP_2(\psi)=\frac {n}2.$$
	Now, for the form $\varphi:=\Pfister{t_{n+1}}\otimes\psi\in I^3F$, which is of dimension 
	$$2(n+2)\leq2\left(\frac d2-2+2\right)=d,$$ 
	we have
	$$GP_3(\varphi)=\frac n2=\left\lfloor \frac d4\right\rfloor-1.$$
	by \ref{PropPfisterZahlUnimodularMalUniformizer} and the conclusion follows.
\end{Beweis}

Furthermore we are already in a good position to determine an upper bound for the 3-Pfister number over rigid fields that generalizes \cite[Theorem 1.13]{Raczek2013}. Our main ingredient is \ref{KorCDVAufRigid}, which was proved with valuation theory.

\begin{Theorem}\label{Schranke3Pfisterzahl}
	Let $F$ be a rigid field. For all even $d\in\N_0$, we have 
	$${GP}_3(F,d)\leq \frac{d^2}{16}.$$
	If we further have $d\geq16$, we even have 
	$${GP_3(F,d)\leq\frac{d^2}{16}-\frac d2-\frac{82-2\cdot(-1)^{\frac d2}}{16}}.$$
\end{Theorem}
\begin{Beweis}
	We will implicitly use that the functions $d\mapsto GP_3(d)$ and $d\mapsto\frac{d^2}{16}-\frac d2-\frac{82-2\cdot(-1)^{\frac d2}}{16}$ are monotonically increasing on the set of even integers $\geq 16$ without referring to this fact explicitly. 
	We use induction on $d$. We already know the following inequalities
	\begin{align*}
		GP_3(F,d)&=0\text{ for all even }d<8,~~
		{GP}_3(F,8)=GP_3(F,10)=1,\\
		{GP}_3(F,12)&=2,~~
		{GP}_3(F,14)=2,  ~~
		{GP}_3(F,16)=3,
	\end{align*}
	that are all compatible with the assertion. 
	As we obviously have the inequality
	$$\frac{d^2}{16}-\frac d2-\frac{82-2\cdot(-1)^{\frac d2}}{16}\leq\frac{d^2}{16}$$
	 for $d\geq16$, we only have to show the second bound.\\
	If a form $\varphi\in I^3F$ of dimension $d\geq 16$ is similar to $d\times \qf 1$ it is Witt equivalent (in fact even isometric) to a sum of $\frac d8$ elements in $GP_3F$ and we are done. 
	Otherwise we can bound ${GP}_3(\varphi)$ according to \ref{KorCDVAufRigid} by
	\begin{align*}
		{GP}_3(F,d-2)+{GP}_2(F,k),
	\end{align*}
	where $k$ is the biggest integer $\leq\frac d2$ that is divisible by two, i.e. we have $$k=\frac d2-\frac12+(-1)^{\frac d2}\cdot\frac 12=2\cdot\left\lfloor\frac d4\right\rfloor,$$
	as we can assume the form $\tau$ in \ref{KorCDVAufRigid} to be of dimension at most $\leq\frac d2$ after possibly scaling with a uniformizer (note that $\tau$ is the second residue class form). By \ref{2PfisterZahl} we thus know
	\begin{align*}
		{GP}_2(F, k)&\leq{GP}_2\left(F, \frac d2-\frac12+(-1)^{\frac d2}\cdot\frac 12\right)\\
		&=\frac{\frac d2-\frac12+(-1)^{\frac d2}\cdot\frac 12}2-1=\frac d4-\frac 54+(-1)^{\frac d2}\cdot \frac14,
	\end{align*}
	which leads to 
	\begin{align}\label{EinmaaaaaaaalUndNieWieeeeeder}
		{GP}_3(F,d)\leq {GP}_3(F, d-2)+\frac d4-\frac 54+(-1)^{\frac d2}\cdot \frac14.
	\end{align}
	We now put $n:=\frac d2-8$, which is equivalent to $d=2n+16$, and consider for $n\in\N$ the recurrence relation
	$$	a_n=a_{n-1}+\frac n2+\frac{11}4+(-1)^n\cdot\frac14,	$$
	which was build by replacing the inequality with an equality in \eqref{EinmaaaaaaaalUndNieWieeeeeder}. For $a_0=3$ (corresponding to $GP_3(F,16) = 3$) this relation has the unique solution
	$$	a_n=\frac18\big(	2n(n+12)+(-1)^n+23	\big)=\frac{d^2}{16}-\frac d2-\frac{82-2\cdot(-1)^{\frac d2}}{16}.	$$
	By construction this is an upper bound for ${GP}_3(\varphi)$ and the proof is complete.
\end{Beweis}

\begin{Bemerkung}
	For non-rigid fields, the 3-Pfister number of quadratic forms may grow exponentially in terms of the dimension, see \cite[Theorem 1.1]{BrosnanReichsteinVistoli} (with \ref{PropVergleichScaledUnscaled} in mind).
\end{Bemerkung}

We can use the above result with an induction to also get upper bounds for the $n$-Pfister numbers of forms in $I^nF$ for any $n\geq4$. We will estimate a little bit coarser to get more succinct bounds. We will further use the following number theoretic result due to Jacob I. Bernoulli \cite{Bernoulli}.

\begin{Theorem}{\cite[Chapter 15, Theorem 1]{IrelandRosen}}
	Let $m\in\N$ be an integer. Then there is some polynomial $p\in\Q[X]$ of degree ${\deg(p)=m+1}$ such that 
	$$	1^m+2^m+\ldots+ n^m=p(n)	$$
	for all $n\in\N$.
\end{Theorem}

Using the distributive rule and the above result several times, we immediately get the following consequence:

\begin{Korollar}\label{IchKannEsNichtAndern}
	Let $q\in\Q[X]$ be a polynomial of degree $\deg(q)=m$. Then there is some polynomial $p\in\Q[X]$ of degree $m+1$ such that we have
	$$	q(1)+q(2)+\ldots +q(n)=p(n)	$$
	for all $n\in\N$.
\end{Korollar}

The main result of this chapter is the following which states that Pfister numbers over all rigid fields can only increase polynomially. For non-rigid fields, it is not even known if the Pfister numbers are finite, see \cite[Remark 4.3]{BrosnanReichsteinVistoli}.

\begin{Theorem}
	Let $n\geq 3$ be an integer. Then there is some polynomial $p\in\Q[X]$ of degree $n-1$ whose associated function $\R_{\geq0}\to\R$ is increasing, nonnegative and fulfils  
	$$GP_n(d, F)\leq p(d)$$ 
	for all rigid fields $F$ and all even integers $d\geq 2^n$.
\end{Theorem}
\begin{Beweis}
	We prove this by induction on $n$, where the induction base $n=3$ is covered by \ref{Schranke3Pfisterzahl}. So let now $n\geq 4$ and let $q_{n-1}\in\Q[X]$ be the polynomial as described in the statement for $n-1$ that exists due to the induction hypothesis and let $p_{n-1}\in\Q[X]$ be the polynomial of degree $n-1$ with 
	\begin{align}\label{MikroPlastikVonSchuhSohlen}
		q_{n-1}(1)+\ldots q_{n-1}(k)={p_{n-1}}(k)
	\end{align}
	for all $k\in\N$ that exists by \ref{IchKannEsNichtAndern}. Obviously, the function $$\R_{\geq0}\to\R, x\mapsto p_{n-1}(x)$$
	is increasing and nonnegative as the function defined by $q_{n-1}$ is so.
	
	Just as in the proof of \ref{Schranke3Pfisterzahl} we have 
	$$GP_n(d,F)\leq GP_n(d-2,F)+GP_{n-1}\left(\frac d2-\frac12+\frac12\cdot(-1)^{\frac d2},F\right)$$ which is - using the same argument again - lower than or equal to
	$$	GP_n(d-4,F)+GP_{n-1}\left(\frac {d-2}2-\frac12+\frac12\cdot(-1)^{\frac {d-2}2},F\right)+GP_{n-1}\left(\frac d2-\frac12+\frac12\cdot(-1)^{\frac d2},F\right).	$$
	Iterating this process, we get a sum of expressions of the form $GP_{n-1}(k, F)$ with $2^{n-1}\leq k\leq \frac d2$ - each of these summands occuring at most 2 times - and one summand of the form $GP_n(2^n, F)$.\\
	As we have $GP_n(2^n, F)=1$ according to the Arason-Pfister Hauptsatz, we thus get the upper bound
	\begin{align*}
		1+2\sum\limits_{k=2^{n-1}, 2\mid k}^{\lfloor\frac d2\rfloor}GP_{n-1}(k,F)&\leq 
		1+2\sum\limits_{k=2^{n-1}, 2\mid k}^{\lfloor\frac d2\rfloor}q_{n-1}(k)\\
		&\leq 1+2\sum\limits_{k=1}^{\lfloor\frac d2\rfloor}q_{n-1}(k)\\
		&=1+2{p_{n-1}}\left(\lfloor\frac d2\rfloor\right)\leq 1+2p_{n-1}\left(\frac d2\right).
	\end{align*}
	It is thus easy to see that the polynomial $p_n(X):=1+2p_{n-1}\left(\frac X2\right)$ does the job.

\end{Beweis}

\section*{Acknowledgments and Notes}
The results contained in this paper are part of the PhD-Thesis of the author. He would like to thank Detlev Hoffmann for supervising this work and giving some very useful hints and corrections in the process.

\bibliographystyle{abbrv}
\bibliography{literatur}

\end{document}